\documentclass[]{interact}

\usepackage{epstopdf}
\usepackage{subfigure}

\usepackage[numbers,sort&compress]{natbib}
\bibpunct[, ]{[}{]}{,}{n}{,}{,}
\renewcommand\bibfont{\fontsize{10}{12}\selectfont}

\theoremstyle{plain}
\newtheorem{theorem}{Theorem}[section]

\newtheorem{conjecture}[theorem]{Conjecture}

\theoremstyle{definition}

\theoremstyle{remark}
\newtheorem{remark}[theorem]{Remark}

\newcommand{\mybrace}[1]{\left\{#1\right\}}
\newcommand{\myspan}[1]{\ensuremath{\operatorname{span}\!\mybrace{#1}}}
\newcommand{\matlab}{MATLAB\textsuperscript{\textregistered}}

\usepackage{url}
\usepackage{graphicx}      
\usepackage{mathtools}	
\usepackage{silence}
\allowdisplaybreaks

\begin{document}

\title{A Connection Between Time Domain Model Order Reduction and Moment Matching for LTI Systems}

\author{
\name{Manuela Hund\textsuperscript{a} and Jens Saak\textsuperscript{a}}
\affil{\textsuperscript{a}Max Planck Institute for Dynamics of Complex
Technical Systems, Sandtorstra{\ss}e 1, 39106 Magdeburg}
}

\maketitle

\begin{abstract}
  We investigate the time domain model order reduction (MOR) framework using
  general orthogonal polynomials by Jiang and Chen~\cite{morJiaC12} and extend
  their idea by exploiting the structure of the corresponding linear system of
  equations. Identifying an equivalent Sylvester equation, we show a connection
  to a rational Krylov subspace, and thus to moment matching. This theoretical
  link between the MOR techniques is illustrated by three numerical examples. For
  linear time-invariant systems, the link also motivates that
  the time domain approach can be at best as accurate as moment matching,
  since the expansion points are fixed by the choice of the polynomial basis, 
  while in moment matching they can be adapted to the system.
\end{abstract}


\section{Introduction}\label{sec:introduction}
Various mathematical and physical processes can be modeled
as linear time-invariant (LTI) input-output systems
\begin{align}
\label{eq:LTI}
   \begin{split}
      E \dot{x}(t) &= A x(t) + B u(t), \\
      y(t) &= C x(t),
   \end{split}
\end{align}
where $E , A \in \mathbb{R}^{n \times n}$ are sparse matrices, $B \in
\mathbb{R}^{n \times p}$ and $C \in \mathbb{R}^{q \times n}$ are input
and output matrices, respectively, $x(t) \in \mathbb{R}^n$ is the state
vector, $u(t) \in \mathbb{R}^p$ is the input vector, $y(t) \in
\mathbb{R}^q$ is the output vector and $t \in \mathbb{R}$ represents
time. 

Since the order of the LTI system~\eqref{eq:LTI} is often huge 
$n \gg 10^3$, a numerical simulation might be too expensive or even 
impossible, caused by immense computational time and memory 
requirements. 
Nevertheless, the input-output behavior of \eqref{eq:LTI} can be computed fast 
and accurately if the given system is reduced to a system
\begin{align}
\label{eq:rLTI}
  \begin{split}
      E_r \dot{x}_r(t) &= A_r x_r(t) + B_r u(t), \\
      y_r(t) &= C_r x_r(t),
   \end{split}
\end{align}
that approximates the dynamic behavior of \eqref{eq:LTI},
but $E_r , A_r \in \mathbb{R}^{r \times r}$, $B_r \in
\mathbb{R}^{r \times p}$, $C_r \in \mathbb{R}^{q \times r}$,
$x_r(t) \in \mathbb{R}^r$, $y_r(t) \in
\mathbb{R}^q$ with the reduced order $r \ll n$. 

The aim of MOR is to approximate a system \eqref{eq:LTI} with a huge order $n$
by a system \eqref{eq:rLTI} with a much smaller order $r$, such that structural 
properties are preserved and the approximation error $y(t)-y_r(t)$ over
a given time interval $[t_{0},t_{f}]$, or the error of the transfer functions
$G(s)-G_r(s)$ over a frequency range $[s_{0},s_{1}]$, is small. 

Transfer functions describe the relation between input and output in frequency 
domain.  
For \eqref{eq:LTI}, \eqref{eq:rLTI} with zero initial states and an
evaluation point $s$ in frequency domain, these are given by
\begin{align}
  G(s) &= C(sE-A)^{-1}B, \label{eq:transfer} \\
  G_r(s) &= C_r(sE_r-A_r)^{-1}B_r. \label{eq:transfer_r}
\end{align}
There exist numerous MOR approaches.
$\mathcal{H}_2$ optimal MOR techniques like the iterative rational
Krylov algorithm (IRKA)  (see, e.g. \cite{morGugAB08}) 
or the two-sided iteration algorithm (TSIA) (see, e.g. \cite{morXuZ11})
measure their approximation error $\epsilon_2$ in the $\mathcal{H}_2$ system norm 
(see, e.g. \cite[Chapter~5]{morAnt05})
\begin{align*}
  \epsilon_2 \coloneqq \Vert G(s)-G_r(s) \Vert_{\mathcal{H}_2}. 
\end{align*}
These techniques are just two types of projection based MOR methods. 
There, the system 
\eqref{eq:LTI} is reduced using so-called projection matrices 
$V,W \in \mathbb{R}^{n \times r}$ that map the matrices $E,A,B$ and 
$C$ onto a subspace approximating the space of the state vector 
$x(t)$. The reduced system is given as
\begin{align}
\label{eq:projection}
  \begin{split}
      \underbrace{W^T E V}_{E_r} \dot{x}_r(t) &= 
      \underbrace{W^T A V}_{A_r}  x_r(t) 
      + \underbrace{W^T B}_{B_r} u(t),\\
      y_r(t) &= \underbrace{CV}_{C_r}  x_r(t),
  \end{split}
\end{align}
where $x(t) \approx V x_r(t)$. A more simple example for projection methods
is moment matching, where $V$ and $W$ are computed to approximate the
moments of the transfer function. This method and its relevant properties
are repeated in Section~\ref{subsec:moment_matching}. It represents one step in
the procedure behind the IRKA iteration.  

Also the balanced truncation technique falls into the class of projection based
reduction methods. Its error $\epsilon_{\infty}$ is measured in the 
$\mathcal{H}_{\infty}$ norm (see, e.g. \cite[Chapter~5]{morAnt05})
\begin{align*}
 \epsilon_{\infty} \coloneqq \Vert G(s)-G_r(s) \Vert_{\mathcal{H}_{\infty}}. 
\end{align*}
Applying this method, the system \eqref{eq:LTI} is first
balanced, i.e. the observability and controllability Gramians
$P_O$ and $P_C$, given as the solutions of two
Lyapunov equations
\begin{align*}
 A^T P_O E   + E^T P_O A   &= - C^TC, \\
 A   P_C E^T + E  P_C A^T &= - BB^T,
\end{align*}
are made equal and diagonal, such that
$P_O=P_C= \text{diag}(\sigma_1 \cdots \sigma_n)$ and
$\sigma_1 \geq \cdots \geq \sigma_n > 0$ are the systems invariant 
Hankel singular values~(HSVs). 
The discardable portions are identified and truncated according to the
magnitude of the HSVs. More details about this method can be found, e.g.
in \cite[Chapter~$7$]{morAnt05}.

The above MOR techniques are motivated and derived by frequency domain
considerations. 

In contrast to that, we next review the idea of Jiang and
Chen~\cite{morJiaC12} presenting a time domain MOR framework
based on orthogonal polynomials. In this paper we
only consider single-input single-output (SISO)
systems, i.e. $p=q=1$ in \eqref{eq:LTI} to simplify the notation.  
Drawbacks of this method are the dependence of the reduced order model (ROM) on
the input $u(t)$ and the initial state $x(t_0)=x_0$.

The dependence on the input can be neglected, since we will see that
piecewise constant controls, which are the most important ones in practical
applications anyway, allow for a joint ROM to exist. However, a ROM
depending on the initial state is undesirable since the reduced system needs to
be recomputed for each initial state or ROMs need to be stored for all
possible initial values. 

Frequency domain based model reduction methods, such as balanced 
truncation or moment matching, assume $x(t_0)=\left[0\ \cdots\  0\right]^T=:~
\mathbb{O}_{n,1}\in \mathbb{R}^n $, in
the first place, in order to avoid additional terms in the transfer function
representation~\eqref{eq:transfer}. For the comparison, we will do the same in
the time domain case in the following. 

The remainder of this paper is organized as follows. The time domain MOR
approach based on general orthogonal polynomials by Jiang and
Chen~\cite{morJiaC12} is introduced in Section~\ref{sec:time_domain}.  In
Section~\ref{sec:exploiting}, the structure of the resulting linear system of
equations is exploited to derive an equivalent Sylvester equation. Further a
slight variation of the approach is discussed that eliminates the initial
condition in the case it is assumed to be zero, which also simplifies the
structure of the coefficients in the Sylvester equation. Since we want to show a
connection to moment matching, we briefly introduce this Krylov subspace
method in Section~\ref{sec:moment_matching} concluding with an important
equivalence to the approaches of Section~\ref{sec:exploiting}. 
Still, the additional freedom in the choice
of the coefficients in the Sylvester equation makes moment matching
theoretically more flexible and better adaptable to the original system
under investigation. Numerical examples illustrated in
Section~\ref{sec:example} demonstrate this advantage. Concluding remarks are
given in Section~\ref{sec:conclusion}.


\section{Time Domain MOR Based on Orthogonal Polynomials (TDMOR)} 
\label{sec:time_domain}
As already mentioned above, we restrict ourselves to SISO systems, i.e. $p=q=1$.

The framework in \cite{morJiaC12} uses  $W=V$ in \eqref{eq:projection} and
obtains the projection matrix $V \in \mathbb{R}^{n \times r}$ from the
vector valued coefficients in series expansions of the state and input, sampling
their time dependence via orthogonal polynomials 
\cite[Chapter~22]{AbrS84}.

The key property of orthogonal polynomials for the derivation of the framework
in \cite{morJiaC12} is given by the following theorem.

\begin{theorem}[Differential recurrence formula, e.g. 
\protect{\cite[Section~$2.15$]{DatM95}}]
\label{thm:diff_rec_formula}
For three subsequent orthogonal polynomials $g_i(t)$ ($i \in \mathbb{N}_0$) holds
\begin{align*}
  g_n(t) &= \alpha_n \dot{g}_{n+1}(t) + \beta_n \dot{g}_n(t) +
  \gamma_n \dot{g}_{n-1} (t), \quad \forall n \in \mathbb{N},
\end{align*}
where $\alpha_n, \beta_n$ and $\gamma_n$ are differential recurrence
coefficients. A list of such coefficients for selected families can be found in
Table~\ref{tab:rec}.
\end{theorem}

\begin{table}
    \tbl{Differential recurrence coefficients}
    {\begin{tabular}{|l||c|c|c|}
      \hline 
      \rule{0pt}{12pt}polynomial class & $\alpha_{i}$ & $\beta_{i}$ &$\gamma_{i}$\\
      \hline 
      \rule{0pt}{12pt}Chebychev-1&$\frac{1}{2i+2}$&0&$-\frac{1}{2i-2}$\\[2ex]
      Chebychev-2&$\frac{1}{2i+2}$&0&$-\frac{1}{2i+2}$\\[2ex]
      Hermite&$\frac{1}{2i+2}$&0&0\\[2ex]
      Jacobi $(a, b>-1)$&
        $\frac{2(a+b+i+1)}{(a+b+2i+2)(a+b+2i+1)}$&
        $\frac{2(a-b)}{(a+b+2i)(a+b+2i+2)}$&
        $-\frac{2(a+i)(b+i)}{(a+b+2i+1)(a+b+2i)(a+b+i)}$\\[2ex]
      Laguerre&-1&1&0\\[2ex]
      Legendre&$\frac{1}{2i+1}$&0&$ -\frac{1}{2i+1}$\\
      \hline
    \end{tabular}}
    \label{tab:rec}
\end{table}

We restrict our considerations to the polynomials investigated in
\cite{morJiaC12}. Other polynomials fulfilling Theorem~\ref{thm:diff_rec_formula} 
are for instance the Gegenbauer polynomials, a generalization of the
Legendre polynomials (see, e.g. \protect{\cite[Section~$2.11$]{DatM95}}.

The following repeats some of the details of the derivation in \cite{morJiaC12}. 
First, the state, the initial condition and the input vector are approximated by 
the following truncated series expansions:
\begin{alignat}{3}
  x(t)         &\approx& x_r(t)  &=& \sum_{i=0}^{r-1} v_i g_i(t),
  \label{eq:approx_x} \\
  x_0 = x(t_0) &\approx& x_r(t_0)&=& \sum_{i=0}^{r-1} v_i g_i(t_0), 
  \label{eq:approx_x0}\\
  u(t)         &\approx& u_r(t)  &=& \sum_{i=1}^{r-1} w_i \dot{g}_i(t) ,
  \label{eq:approx_u}
\end{alignat}
where $v_i \in \mathbb{R}^n$ and $w_i \in \mathbb{R}$ are weights 
determining the subspace $\myspan{V}$ and $g_i(t)$ are orthogonal 
polynomials representing the time dependence.
Note, that the open literature provides no information about the remainder terms
in \eqref{eq:approx_x}-\eqref{eq:approx_u}. Thus, the estimation of
approximation errors, and resulting model reduction errors, is at best difficult.

The approximations of the state \eqref{eq:approx_x} and the 
input \eqref{eq:approx_u} are inserted into the state equation of \eqref{eq:LTI}
and using Theorem~\ref{thm:diff_rec_formula}, one 
obtains an expression that only depends on $\dot{g}_i(t)$, since 
$g_0$ is always constant:
\begin{align*}
    B \left( \sum_{i=1}^{r-1} w_i \dot{g}_i(t) \right) =& \sum_{i=1}^{r-1} 
    \left( E - \beta_i A \right) v_i \dot{g}_i(t) -Av_0 g_0(t) - 
     \sum_{i=2}^{r} \alpha_{i-1} A v_{i-1} \dot{g}_i(t) - \\
    & -\sum_{i=1}^{r-2} \gamma_{i+1} A v_{i+1} \dot{g}_i(t).
\end{align*}
A comparison of coefficients leads to the huge $(nr \times nr)$ linear system of
equations $H v= f$ presented in equation \eqref{eq:GLS}, where the approximation
of the initial state \eqref{eq:approx_x0} is only required to obtain a square
matrix

\begingroup\makeatletter\def\f@size{8}\check@mathfonts
\def\maketag@@@#1{\hbox{\m@th\large\normalfont#1}}%
\begin{align}
\label{eq:GLS}
    \begin{bmatrix}
       g_0(t_0)I_n                   & g_1(t_0)I_n   & g_2(t_0)I_n    & g_3(t_0)I_n     & \cdots         & g_{r-1}(t_0)I_n     \\     
      - \frac{g_0(t)}{\dot{g}_1(t)}A& E - \beta_1 A & - \gamma_2 A   &                 &                &          \\
                                    & -\alpha_1 A   & E -\beta_2 A   & -\gamma_3 A     &                &          \\
                                    &		 & \ddots         & \ddots          & \ddots         &           \\
                                    &               &                & -\alpha_{r-3}A  & E-\beta_{r-2}A & -\gamma_{r-1} A  \\
                                    &               &                &                 & -\alpha_{r-2}A & E-\beta_{r-1} A 
    \end{bmatrix}
    \begin{bmatrix}
      v_0 \vphantom{g_0(t_0)I_n}\\
      v_1 \vphantom{\frac{g_0(t)}{\dot{g}_1(t)}A E - \beta_1 A } \\ 
      v_2 \vphantom{\gamma_2}\\
      \vdots   \vphantom{\ddots}\\
      v_{r-2} \vphantom{\beta_{r-2}} \\ 
      v_{r-1} \vphantom{\beta_{r-1}}
    \end{bmatrix}
    =
    \begin{bmatrix}
      \vphantom{g_0(t_0)I_n} x_0 \\
      \vphantom{\frac{g_0(t)}{\dot{g}_1(t)}A E - \beta_1 A } B w_1 \\ 
      \vphantom{\gamma_2} B w_2 \\
      \vdots \vphantom{\ddots}\\
      \vphantom{\beta_{r-2}} B w_{r-2}\\ 
      \vphantom{\beta_{r-1}} B w_{r-1}
    \end{bmatrix},
\end{align}
\endgroup
and $I_n$ denotes the $n$-dimensional identity matrix.

In \cite{morJiaC12} this linear system of equations, with matrix
$H \in \mathbb{R}^{nr \times nr}$ and right hand side $f \in \mathbb{R}^{nr}$, 
is solved using an iterative algorithm. The solution vector 
$v \in \mathbb{R}^{nr}$ is then used to compute the projection matrix 
$V \in \mathbb{R}^{n \times r}$ by orthogonalizing the span of
$\left[ v_1 , \cdots , v_r \right]$. 
In the following context, we will call this method
TDMOR.

Note that $H$ is not depending on time,
since both $g_{0}$ and $\dot g_{1}$ are constant in time.
Further, the matrix $H$ has a certain block-structure. We exploit this structure
in the
following section to derive an equivalent formulation and a more well-posed variation of this
MOR method.  


\section{Structure Exploitation and a Slight Variation} 
\label{sec:exploiting}

\subsection{Structure Exploitation (\emph{SYLTDMOR1})}\label{subsec:original}
 
We multiply the first equation in \eqref{eq:GLS} by $A$
and obtain the following equivalent linear system of equations in Kronecker
product (see, e.g. \protect{\cite[Section~$4.2$]{HorJ91}}) form
\begin{align*}
 \left( \tilde{E}^T \otimes A + \tilde{A}^T \otimes E
 \right)v=\tilde{f}.
\end{align*}
Here 
\begin{align*}
\tilde{E}^T &=
\begin{bmatrix} 
 g_0(t_0)                     & g_1(t_0)      & g_2(t_0)       & g_3(t_0)     & \cdots             & g_{r-1}(t_0)    \\     
- \displaystyle \frac{g_0(t)}{\dot{g}_1(t)} & -\beta_1       & - \gamma_2    &           &          &       \vphantom{\vdots}             \\
     & -\alpha_1     & -\beta_2       & -\gamma_3    &         &                  \\
     &         & \ddots   & \ddots & \ddots         &                 \\
     &          &        & -\alpha_{r-3}  & -\beta_{r-2}   & -\gamma_{r-1}   \\
     &          &        &                & -\alpha_{r-2}  & -\beta_{r-1} \vphantom{\vdots}
 \end{bmatrix} \in \mathbb{R}^{r \times r},\\
  \tilde{A}^T &= 
    \begin{bmatrix}
      \mathbb{O}_{1,1} & \mathbb{O}_{1,r-1}  \\
      \mathbb{O}_{r-1,1} & I_{r-1} \\
    \end{bmatrix} \in \mathbb{R}^{r \times r} , \\
  \tilde{f}&=
    \begin{bmatrix}
      (A x_0)^T & (B w_1)^T & \hdots & (B w_{r-1})^T
    \end{bmatrix}^T \in \mathbb{R}^{n r}.
\end{align*}

Using the equivalence 
\begin{align}
\label{eq:equiv_GLS}
 \left(B^T \otimes A \right) \text{vec}(X)= \text{vec}(C) 
 \Leftrightarrow AXB=C,
\end{align}
 (see, e.g. in \cite[Section~$4.3$]{HorJ91}) we obtain a 
 Sylvester equation
\begin{align}
\label{eq:Sylv}
 AV\tilde{E} + EV\tilde{A} = \tilde{F} ,
\end{align}
where $v= \text{vec}(V) $, $\tilde{f}= \text{vec}\left(\tilde{F}\right) $
and vec$(.)$ of a matrix $A \in \mathbb{R}^{m \times n}$ is defined as in
\begin{align*}
   \text{vec} (A) =
  \begin{bmatrix}
    a_{1,1}, \hdots, a_{m,1}, a_{1,2}, \hdots, a_{m,2}, 
    \hdots, a_{1,n},\hdots, a_{m,n}
  \end{bmatrix}^T
\end{align*}

\noindent and is called a \textit{vectorization}
(see, e.g. \protect{\cite[Section~$4.2$]{HorJ91}}).
A further orthogonalization of $V$ leads to the desired
projection matrix. In the following, we will refer to this algorithm as
SYLTDMOR1.

Note that the matrix pencil $(\tilde{A}, \tilde{E})$ has at least one eigenvalue
equal to zero caused by the structure of $\tilde{A}$, which arises from the 
initial state condition. 
Therefore the matrix $H$ in \eqref{eq:GLS} is not invertible resulting in an 
infinite number of solutions and thus in a possibly infinite number of ROMs. 
Hence, this method is not well-posed and its solution not well-defined.


\subsection{Variation of the Presented Algorithm (\emph{SYLTDMOR2})}
\label{subsec:variation}
In the approach of Jiang and Chen \cite{morJiaC12}, the approximation of the
initial state is only required to obtain a square matrix. 
As a conclusion of Section \ref{subsec:original}, this condition turned out to 
be linearly dependent anyway. Besides, a ROM depending on the initial state is
not desirable.

In Section~\ref{sec:example}, we compare to the frequency domain methods, thus
also here we fix the initial state to $x_0= x(t_0)= \mathbb{O}_{n,1}$.
 
Doing so, we can neglect the constant polynomials $g_0(t)$ in the
approximations \eqref{eq:approx_x}. 
In order to keep an $r$ dimensional approximation, we shift the sums by 1.
Using the same procedure as in Section~\ref{sec:time_domain},
we end up with an ${nr \times nr}$ linear system of
equations $\hat{H}\hat{v}= \hat{f}$. Rewriting it, again using the Kronecker product, we obtain 
\begin{align*}
 \left( \hat{E}^T \otimes A + \hat{A}^T \otimes E \right) \hat{v} = \hat{f},
\end{align*}
where now
\begin{align*}
  \hat{E}^T &= - 
    \begin{bmatrix}  
      \beta_1     &  \gamma_2   &                 &              &     \\
      \alpha_1    & \beta_2     & \gamma_3        &              &     \\
		  & \ddots      & \ddots          & \ddots       &     \\
		  &             & \alpha_{r-2}    & \beta_{r-1}  & \gamma_{r}  \\
		  &             &                 & \alpha_{r-1} & \beta_{r} 
    \end{bmatrix} \in \mathbb{R}^{r \times r} , \qquad
    \hat{A}^T = I_r, \\
    \hat{v} &=
    \begin{bmatrix}
      \hat{v}_1^T & \hdots & \hat{v}_r^T
    \end{bmatrix} \in \mathbb{R}^{nr}, \qquad
      \hat{f} =
    \begin{bmatrix}
      \left( B w_1 \right)^T & \hdots & \left( B w_r \right)^T
    \end{bmatrix}^T \in \mathbb{R}^{nr}.
\end{align*}

Exploiting the equivalence \eqref{eq:equiv_GLS} and the fact $\hat{A}=I_r$,
the linear system of equations $\hat{H}\hat{v}=\hat{f}$ can be 
reformulated as the Sylvester equation
\begin{align}
\label{eq:Sylv_2}
 A\hat{V} \hat{E} + E \hat{V} = \hat{F},
\end{align}
where $\hat{v}= \text{vec}\left(\hat{V}\right) $ and $\hat{f}= 
\text{vec}\left(\hat{F}\right) $. As in Section~\ref{subsec:original},
the projection matrix can be obtained by orthogonalization of \(\hat V\). 
In the following, we will call this method SYLTDMOR2.

Compared to \eqref{eq:Sylv}, Sylvester equation
\eqref{eq:Sylv_2} does not depend on the initial state.
Moreover, the pencil $(I_r,\hat E)$ does not have a zero eigenvalue, such that (in
contrast to \eqref{eq:Sylv}) \eqref{eq:Sylv_2} always allows for
a unique solution. 
Thus, this method is well-posed and the ROM is well-defined.


\subsection{Reincorporation of Non-Zero Initial Conditions in \emph{SYLTDMOR2}}
\label{sec:inclusion-non-zero}

Although we removed it in the formulation, it is possible to use the initial
state condition in SYLTDMOR2.
One way to include the initial condition is given by the approach presented
in \cite{morHeiRA11}, where the given SISO system is reformulated to a 
multiple-input single-output (MISO) system, by adding the initial state as a
column in \(B\) and using a corresponding Dirac input. 
Another and more flexible method is described in \cite{morBeaGM17} for the 
frequency domain MOR methods, we want to compare with.
Here, a whole variety of initial state conditions, instead of only one 
condition, can be considered using an approach splitting the problem into a
homogeneous and inhomogeneous part, that can be solved separately. 
This method preserves the SISO system and can also be applied to the time 
domain MOR approach.
If the subspace of relevant initial conditions is known, this method clearly
offers a more flexible setting and overcomes the problem of storing a separate
reduced model for every possible initial condition.

\section{Moment Matching and its Relation to \emph{SYLTDMOR2}} 
\label{sec:moment_matching}

Our main goal in
this paper is to show a connection between the above mentioned time domain MOR
approaches and moment matching. To this end, we repeat the basics of this Krylov
subspace technique by first introducing a standard Krylov subspace
(see, e.g. \cite[Section~1.6]{morWol15}) of order $r$ for a matrix $A \in \mathbb{R}^{n 
\times n}$ and a vector $b \in \mathbb{R}^n$ as
 \begin{align*}
 \mathcal{K}_r(A,b)= \myspan{b, Ab, \hdots, A^{r-1}b} .
 \end{align*}

\subsection{Moment Matching} \label{subsec:moment_matching} 
Moment matching is a projection based MOR technique. It constructs the
projection matrix starting from a series expansion of the transfer 
function rather than the state exploiting the Neumann series (see, e.g. \cite{Yos95})
\begin{align}
\label{eq:Neumann}
  \left(I_n-T\right)^{-1}=\sum_{k=0}^{\infty} T^k,
\end{align} 
where $T\in \mathbb{R}^{n \times n}$ is a matrix, such that 
$(I_n-T)$ is in fact invertible.

Assuming, that $(s_0 E -A)$ is invertible, and using
\eqref{eq:Neumann}, the transfer function of the original system
can be expressed as
\begin{align*}
    G(s) =   \sum\limits_{k=0}^{\infty}
    \underbrace{ 
           C \left( -(s_0 E-A)^{-1}E  \right)^k (s_0 E-A)^{-1} 
           B}_{M_k^{s_0}}  (s-s_0)^k,
\end{align*}
where $M_k^{s_0}$ are called moments of the original transfer function around
$s_0$.

The aim of moment matching is to find a reduced system of 
order $r \ll n$, such that for some $k= 0, \hdots ,\infty$ for the moments of the reduced
order transfer function we have $\hat{M}_k^{s_0}=M_k^{s_0}$ .

This equality of moments can be guaranteed by using an orthonormal basis of the
input or output Krylov subspace around a single expansion point $s_0 \in
\mathbb{C}$ to form the orthogonal matrices $Q_1$ and $Q_2$
\begin{align*}
  \mathcal{K}_m \left((A-s_0 E)^{-1}E, (A-s_0 E)^{-1}B \right) &=
                                           \myspan{Q_1},\\
  \mathcal{K}_m \left((A-s_0 E)^{-T}E^T, (A-s_0 E)^{-T}C^T \right) &=
                                                              \myspan{Q_2}.
\end{align*}
If the one-sided Krylov subspace method is used, i.e. $V=W=Q_1$ is used to
project, $r$ moments will match (see, e.g. \cite[Chapter~11]{morAnt05}). In
\cite[Chapter~3]{morEid09} it is pointed out, that this property also holds
if $Q_2$ is used instead of $Q_1$. In contrast, if both $V=Q_1$ and $W=Q_2$,
then $2r$ moments of the original and reduced order systems will match
(see, e.g. \cite[Chapter~11]{morAnt05}). This method is called two-sided Krylov
subspace method.

If multiple expansion points $s_1, \hdots , s_k \in \mathbb{C}$ are given, 
$Q_1$ and $Q_2$ can be obtained as a basis of the union of
Krylov subspaces, that belong to the expansion points:
\begin{align*}
  \bigcup\limits_{i=1}^k \mathcal{K}_{r_i} \left((A-s_i E)^{-1}E, 
  (A-s_i E)^{-1}B \right) &= \myspan{Q_1}, \\
  \bigcup\limits_{i=1}^k \mathcal{K}_{r_i} \left((A-s_i E)^{-T}E^T, 
  (A-s_i E)^{-T}C^T \right) &= \myspan{Q_2},
\end{align*}
where $\sum\limits_{i=1}^k r_i=r$. Using only $V=W=Q_1$ to project,
the first $r_i$ moments around $s_i$ of the original and reduced order
model match for $i=1, \hdots, k$. In the two-sided Krylov subspace method
using multiple expansion points $2 r_i$ moments will match around $s_i$ 
for $i=1,\hdots,k$ (see, e.g. \cite[Chapter~3]{morEid09}).


\subsection{Moment Matching and Sylvester Equations}
Since TDMOR presented in Section~\ref{sec:time_domain} and SYLTDMOR1 and SYLTDMOR2
presented in Section~\ref{sec:exploiting}
only use one projection matrix $V$ to obtain a ROM, we
will only focus on the one-sided Krylov subspace method. 

On the one hand, the projection matrix \(V\) can be obtained using the approach
presented in Section~\ref{subsec:moment_matching}. On the other hand, there is a
very useful result describing a relation between the basis of a Krylov subspace
and the solution of a Sylvester equation, that can be found in
\cite[Section~$3.4$]{morVan04} and \cite[Section~$2.3$]{morWol15}.  This
connection, requires the observability of a matrix pair $(S,L)\in 
\mathbb{C}^{r\times r}\times \mathbb{C}^{p \times r}$. This is, e.g., given (see,
e.g. \cite[Chapter~4]{morAnt05}), when the corresponding
observability matrix
\begin{align*}
  Ob(S,L)=
  \begin{bmatrix}
    L^H & \left(LS \right)^H & \left(LS^2\right)^H & \hdots & \left(LS^{r-1}\right)^H
  \end{bmatrix}^H
\end{align*}
has full rank. One then has the following theorem.

\begin{theorem}[Single expansion point duality, e.g. \protect{\cite[Section~$2.3$]{morWol15}}]\label{thm:equiv_single}
  Given the expansion point $s_0 \in \mathbb{C}$, such that $s_0$ is not 
  an eigenvalue of $E^{-1}A$, the columns of $V \in \mathbb{C}^{n \times r}$
  form a basis of a rational Krylov subspace 
  \begin{align*}
    \myspan{V}= \mathcal{K}_r \left((A-s_0 E)^{-1}E, (A-s_0 E)^{-1}
    B \right),
  \end{align*}
  if and only if there exists an observable pair $(S,L)$, where
  $ S~\in~\mathbb{C}^{r \times r}$, $L~\in~\mathbb{C}^{1 \times r}$, which admits
  the Jordan canonical form $J$,
  \begin{align*}
    T^{-1}ST=J =
    \begin{bmatrix}
      s_0 & 1      &        &    \\
      & \ddots & \ddots &    \\
      &        & \ddots & 1  \\
      &        &        & s_0
    \end{bmatrix},
  \end{align*}
  for an appropriate transformation matrix $T \in \mathbb{C}^{r \times r}$, such
  that the Sylvester equation 
  \begin{align}
   \label{eq:sylv_thm1}
    AV - EVS = BL
  \end{align}
  is satisfied.\\
  Moreover, the reduced model $G_r(s)=C_r (s E_r-A_r)^{-1}B_r$ from 
  \eqref{eq:rLTI} matches the  moments $M_i^{s_0}= \hat{M}_i^{s_0}, 
  i= 0, \hdots , r-1$, if $s_0$ is not a pole of $G_r(s)$.
\end{theorem}

This theorem also extends to the case of multiple expansion points.
\begin{theorem}[Multiple expansion point duality, e.g. 
  \protect{\cite[Section~$2.3$]{morWol15}}]\label{thm:equiv_multiple}
  Given $r$ distinct expansion points $s_1, \hdots, s_r \in \mathbb{C}$, such
  that none of them is an eigenvalue of $E^{-1}A$, the columns of
  $V \in \mathbb{C}^{n \times r}$ form a basis of a rational Krylov subspace
  \begin{align*}
    \myspan{V}=\myspan{(A-s_1E)^{-1}B, \hdots, (A-s_r E)^{-1}B},
  \end{align*}
  if and only if there exists an observable pair $(S,L)$ with
  $ S~\in~\mathbb{C}^{r \times r}$, $L~\in~\mathbb{C}^{1 \times r}$, which
  admits the Jordan canonical form $J$,
  \begin{align*}
    T^{-1}ST=J= \textnormal{diag}(s_1, \hdots, s_r) \textnormal{ and } LT=
    \begin{bmatrix}
      1 & \hdots & 1
    \end{bmatrix}
  \end{align*}
  for an appropriate transformation matrix $T\in \mathbb{C}^{r \times r}$,
  such that the Sylvester equation 
  \begin{align}
  \label{eq:sylv_thm2}
    AV - EVS= BL
  \end{align}
  is satisfied. 
  
  Moreover, the reduced model $G_r(s)=C_r (s E_r-A_r)^{-1}B_r$ from 
  \eqref{eq:rLTI} matches the  moments $M_0^{s_i}= \hat{M}_0^{s_i}, 
  i= 0, \hdots, r-1$, if none of the $s_i$ is a pole of $G_r(s)$.
\end{theorem}

Theorems~\ref{thm:equiv_single} and \ref{thm:equiv_multiple}
describe an important connection between 
a Krylov subspace MOR technique and the solution of a Sylvester equation.
Every basis of a rational Krylov subspace solves a certain 
Sylvester equation consisting of an observable matrix pair $(S,L)$. Here,
the eigenvalues of $S$ correspond to the expansion points in moment matching.
 Following \cite[Theorem $3.23$]{morVan04}, 
 the eigenvalues of $S$ are either interpolation points between $G(s)$ and $G_r(s)$
 or the inverse of common poles between $G(s)$ and $G_r(s)$.
Considering multiple-input multiple-output (MIMO) systems, the matrix $L$ is of 
importance, since tangential directions are stored in its columns.
Conversely, every solution of a Sylvester equation consisting of an
observable matrix pair $(S,L)$ spans a Krylov subspace with expansion
points given by the eigenvalues of $S$.

The following Section~\ref{subsec:equiv_moment} uses
Theorems~\ref{thm:equiv_single} and \ref{thm:equiv_multiple} to show a novel
connection between moment matching and the time domain MOR framework based on
orthogonal polynomials.


\subsection{Equivalence of  \emph{SYLTDMOR2} and Moment Matching}
\label{subsec:equiv_moment}
In this section, we apply Theorems~\ref{thm:equiv_single} and
\ref{thm:equiv_multiple} to the derived Sylvester equations~\eqref{eq:Sylv} and
\eqref{eq:Sylv_2} from Section~\ref{sec:exploiting}. In the moment matching MOR
method, it is assumed, that the initial state vector is
$x_0= x(t_0)= \mathbb{O}_{n,1}$. Since this condition is also required for 
SYLTDMOR2, we only need to set the initial state vector to zero for the remaining 
time domain MOR approaches to compare these methods. Recall, that the approximation of the 
initial state was only needed to derive a square linear system of equations. For a
consistent initial state it is thus redundant and the matrix $H$ in the linear
system \eqref{eq:GLS} is actually singular. 
Another restriction, we make in this paper, is to set (without
loss of generality) the time interval to $t=[0,1]$. Note that for the general
case $t_1 \in [t_0,t_f]$, this can always be obtained by the simple transformation 
$t_1 \mapsto \displaystyle \frac{t_1-t_0}{t_f-t_0}$ for constant time increments.

To obtain the structure of the Sylvester equations \eqref{eq:sylv_thm1} or 
\eqref{eq:sylv_thm2} from Theorems~\ref{thm:equiv_single} and 
\ref{thm:equiv_multiple}, it is necessary to invert the $\tilde{E}$ (SYLTDMOR1) and 
$\hat{E}$ (SYLTDMOR2) matrices containing information about the orthogonal
polynomials.
Due to the structure of these matrices, it is only possible to invert them
in the following cases:

Matrix $\tilde{E}$ is regular for:
\begin{itemize}
\item Hermite: $r$ odd (otherwise $g_{r-1}(t_0)=0$ and thus we obtain a
  zero row)
\item Laguerre: all $r$
\item Legendre, Chebychev of first and second kind: $r$ odd (otherwise a zero row
  is obtained due to linear combination)
\end{itemize}

Matrix $\hat{E}$ is
\begin{itemize}
\item Hermite: always singular
\item Laguerre: always regular
\item Legendre, Chebychev of first and second kind: regular for $r$ even 
(otherwise a zero row is obtained due to linear combination)
\end{itemize}

Explicit representations of the inverse matrices for the different 
polynomials listed above can be found in~\cite{morHun15}.
The inverse matrices of the Jacobi polynomials cannot be obtained as easy
as for the above mentioned polynomials caused by the structure and the influence
of parameters $a$ and $b$. Therefore we assume to choose $a$ and $b$, such that
$\tilde{E}$ and $\hat{E}$ are invertible. 
In the following, the Jacobi polynomials are only used to proof the assumptions 
of Theorems \ref{thm:equiv_single} and \ref{thm:equiv_multiple} since the 
Legendre and Chebychev polynomials are special cases of these polynomials 
(see, e.g. \cite[Chapter~22]{AbrS84}).

Assuming either of the aforementioned cases and exploiting the zero initial
state,  we rewrite the Sylvester equation \eqref{eq:Sylv} as
\begin{alignat*}{2}
  AV \tilde{E} + EV \tilde{A} & =  
    \begin{bmatrix}
	Ax_0 & Bw_1 & \hdots & Bw_{r-1}
    \end{bmatrix}\\
 \Leftrightarrow \quad AV - EVS &= BL,
\end{alignat*}
where $S= - \tilde{A}\tilde{E}^{-1}$ and $L =
\begin{bmatrix}
  z_0 & w_1 & \hdots & w_{r-1}
\end{bmatrix}
\tilde{E}^{-1} $ and $z_0= 0$, since $Ax_0= \mathbb{O}_{n,1}$ due to 
the initial state.

Equivalently, we can rewrite the Sylvester equation \eqref{eq:Sylv_2}:
\begin{alignat}{2}
  \nonumber 
 A \hat{V} \hat{E} + E\hat{V} &= 
      \begin{bmatrix}
	  Bw_1 & \hdots & Bw_{r}
      \end{bmatrix} \\
  \label{eq:sylv_normal}
  \Leftrightarrow \quad  A\hat{V} - E\hat{V}\hat{S}&= B\hat{L},
\end{alignat}
where $\hat{S}=-\hat{E}^{-1}$ and $\hat{L}= 
\begin{bmatrix}
  w_1 & \hdots & w_{r}
\end{bmatrix}\hat{E}^{-1}$.

Since the eigenvalues of $S$ and $\hat{S}$ are the expansion points only in 
case of observability, we now have to check the
observability of the matrix pairs $(S,L)$ and $(\hat{S}, \hat{L})$.  While $S$
and $\hat{S}$ only depend on the choice of the orthogonal polynomial, $L$ and
$\hat{L}$ additionally depend on the expected input $u(t)$, since
$w_1, \hdots, w_r$ are weights of the approximated input \eqref{eq:approx_u}. 

\begin{remark}
In practice, the input often needs to be realized
piecewise constant. Therefore we assume, that $u(t)=1$. Note that any other 
constant value for \(u(t)\) only scales the solution and thus changes the basis
but not the subspace spanned by \(V\). 
As a consequence the reduced order model stays the same.
\end{remark}

\begin{figure}[t]
  \centering
  \includegraphics{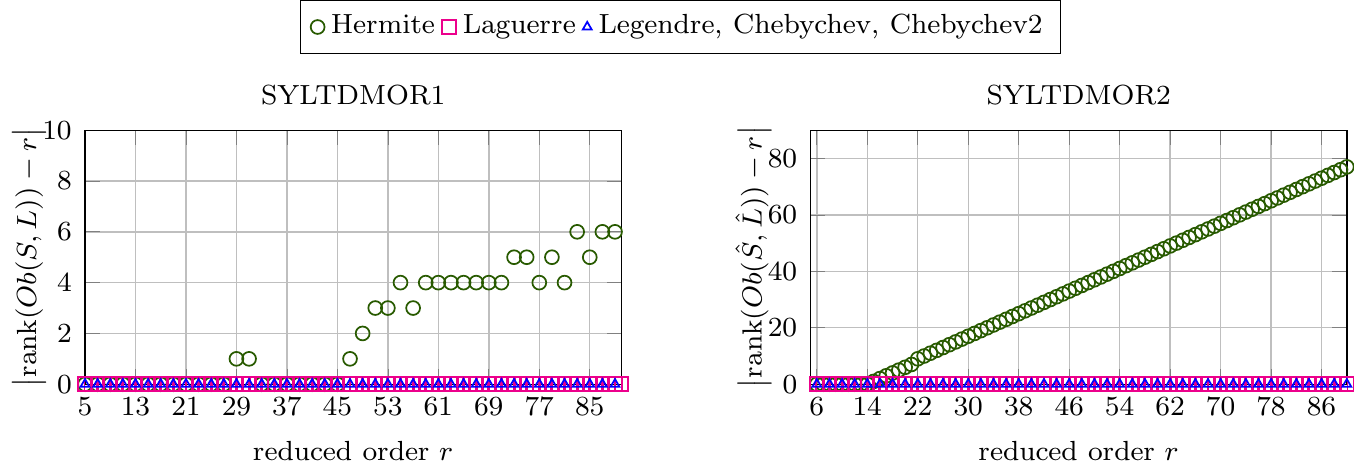}
  \caption{Difference between reduced order and numerical rank of the 
	  observability matrices}
  \label{fig:rank}
\end{figure}

Under this condition we compare the numerical ranks of the associated
observability matrices for certain orthogonal polynomials. Here, the tolerance 
of \matlab s rank function was set to $10^{-20}$ to ensure a good rank 
estimation\footnote{machine precision, i.e. a tolerance \(\approx 10^{-16}\),
  turned out to give unreliable rank decisions in the numerical experiments}. The differences
between the reduced orders $r$ and the numerical rank of the observability
matrix are depicted in Figure~\ref{fig:rank}. 
For the matrix pair $(S,L)$ we only consider an odd order
$r$, for $(\hat{S}, \hat{L})$ only an even order $r$ due to the invertibility
conditions for the matrices $\tilde{E}$ and $\hat{E}$. There are two exceptions:
Since the $\tilde{E}$ and $\hat{E}$ matrices are always invertible in case of
Laguerre polynomials, their numerical rank is plotted for all $r$. 
Even though the $\hat{E}$ matrix using the Hermite polynomials is not invertible
at all, we use here matrix pair $(\tilde{S}, \tilde{L})$ from equation 
\eqref{eq:sylv_new} and thus observability matrix \eqref{eq:obs_hermite} 
instead to obtain its numerical rank.

In both subfigures it is easy to see, that the Legendre and both types of
Chebychev polynomials always lead to a full numerical rank. 
In case of these polynomials,
the rank is only plotted with one mark, because the result is always the same. The
Laguerre polynomials show the same behavior in both figures for all $r$.
In contrast, the Hermite polynomials in SYLTDMOR2 only have a full numerical
rank if the reduced order $r$ is small enough, i.e. $r \leq 14$. For SYLTDMOR1,
their numerical rank is full only if $r \leq 27$.

Note, that Figure~\ref{fig:rank} only presents numerical ranks. 
We now, considering only our
proposed new variant SYLTDMOR2, prove the full rank of the
observability matrices for the Jacobi, Laguerre, Legendre and Chebychev
polynomials of first and second kind rigorously. We will also show, that for the
Hermite polynomials, theoretically, the rank of $Ob(\tilde{S}, \tilde{L})$ is always
full as opposed to the numerical rank.
 
To overcome the difficulties with the singularity of $\hat{E}$ for some kinds
of orthogonal polynomials, we will rewrite Sylvester equation \eqref{eq:Sylv_2}
as 
\begin{align}
\label{eq:sylv_new}
    A \hat{V} \hat{E} + E \hat{V} = B 
    \begin{bmatrix}
      w_1 & \hdots & w_r
    \end{bmatrix}
    \quad \Leftrightarrow \quad   
    E \hat{V} - A \hat{V}\tilde{S} = B \tilde{L},
\end{align}
where $\tilde{S}= - \hat{E} $ and $\tilde{L} =
\begin{bmatrix}
  w_1 & \hdots & w_r
\end{bmatrix}
$, i.e. with the roles of \(E\) and \(A\) swapped.

Let $s_1, \hdots ,s_r \in \mathbb{C} \backslash \Lambda \left(E^{-1}A\right)$
be distinct expansion points, then
\begin{align*}
   \text{span} \big \{ &(A-s_1 E)^{-1}B, \hdots , (A-s_r E)^{-1}B
     \big\} \\ 
   &= \text{span} \Bigg\{- \frac{1}{s_1} \left(E-\frac{1}{s_1}A
     \right)^{-1} B , \hdots , - \frac{1}{s_r}
     \left(E-\frac{1}{s_r}A \right)^{-1} B \Bigg \} \\  
   &= \text{span}\Bigg \{  \left(E-\frac{1}{s_1}A \right)^{-1} B , \hdots
            , \left(E-\frac{1}{s_r}A \right)^{-1} B \Bigg \}.
\end{align*}
Thus, the solution of the Sylvester equation \eqref{eq:sylv_new} is a basis of 
a Krylov subspace with expansion points
$\displaystyle \frac{1}{s_i}$ for $i=1,\hdots,r$.

\subsubsection{Hermite Polynomials}
Since for these polynomials the \(\hat{E} \) matrix is singular for all
$r$, we choose Sylvester equation \eqref{eq:sylv_new} to prove the equivalence
to moment matching. As mentioned above, we assume the input to be chosen
piecewise constant and thus $\tilde{L} =
\begin{bmatrix}
 1 & 0 & \hdots & 0
\end{bmatrix}
$. Further, 
\begin{align*}
 \tilde{S} = - \frac{1}{2} 
 \begin{bmatrix}
  0 & \displaystyle \frac{1}{2}  &        &  \\
    & \ddots                     & \ddots &  \\
    &                            & \ddots & \displaystyle \frac{1}{r}\\
    &                            &        & 0 
 \end{bmatrix}.
\end{align*}

Due to this special structure, the $(r \times r)$ observability matrix is
a diagonal matrix with non-zero entries
\begin{align}
    \label{eq:obs_hermite} 
    \text{diag}(Ob(\tilde{S}, \tilde{L}))= 
    \begin{bmatrix}
	1 & 
	\displaystyle - \frac{1}{2} \cdot \frac{1}{ 2} & 
	\displaystyle   \frac{1}{2^2} \cdot \frac{1}{ 2 \cdot  3} & 
	\hdots &
	\displaystyle \left(-\frac{1}{2}\right)^{r-1} \left( \prod\limits_{i=1}^{r-1} (i+1)\right)^{-1} 
    \end{bmatrix}
\end{align}
and thus $(\tilde{S}, \tilde{L}) $ is observable. 
Converting this problem back to Sylvester equation \eqref{eq:sylv_normal},
the expansion points are generalized eigenvalues of $(-\hat{E}, I_r)$
and thus inverse eigenvalues of $(-I_r, \hat{E})$, where $\hat{E} = - 
\tilde{S}$. Since $\hat{E}$ is a strict upper triangular matrix, all eigenvalues
of $(-I_r, \hat{E})$ are zero and thus all expansion points
are $\infty$ (see, e.g. \cite{morGri97}, \cite[Chapter~4]{Dem97}).

\subsubsection{Laguerre Polynomials}
\label{subsec:Laguerre}
Since the $\hat{E}$ matrix for these polynomials is always invertible,
we choose Sylvester equation \eqref{eq:Sylv_2}.
The explicit inverse to
\begin{align*}
 \hat{E} =
 \begin{bmatrix}
  -1 & 1     &        &    \\
     &\ddots & \ddots &       \\
     &       & \ddots & 1 \\
     &       &        &-1 \\
 \end{bmatrix}
\end{align*}
is given by the upper triangular matrix
\begin{align*}
 \hat{E}^{-1} = -
 \begin{bmatrix}
  1 & \cdots & 1 \\
    & \ddots & \vdots \\
    &        & 1
 \end{bmatrix}.
\end{align*}

Since  $\hat{S}= -\hat{E}^{-1}$ by definition and due to the choice of piecewise 
constant inputs, we have
\begin{align*}
  \hat{L}&=
           \begin{bmatrix}
             1 & 0 & \hdots & 0
           \end{bmatrix} \hat{E}^{-1} 
         = - 
           \begin{bmatrix}
             1 & \hdots & 1
           \end{bmatrix} 
         = - \hat{S}(1,:).
\end{align*}

Thus the entries of the observability matrix become
\begin{align*}
 Ob(\hat{S}, \hat{L})
  =
  \begin{bmatrix}
      \hat{L} \\  \hat{L} \hat{S} \\ \vdots \\  \hat{L} \hat{S}^{r-1}
  \end{bmatrix}
  = -
  \begin{bmatrix}
      \hat{S}(1,:) \\ \hat{S}(1,:)\hat{S} \\ \vdots \\ \hat{S}(1,:)\hat{S}^{r-1}
  \end{bmatrix}
  = -
  \begin{bmatrix}
      \hat{S}(1,:) \\ \hat{S}^2(1,:) \\ \vdots \\ \hat{S}^{r}(1,:)
  \end{bmatrix}.
\end{align*}
Observing that the first rows in the powers of \(\hat{S} \) can be written in terms
of binomial coefficients and using the sum formula 
\begin{align*}
 \sum_{k=0}^n \binom{k}{l} = \binom{n+1}{l+1},
\end{align*}
for integers $k,l,n \geq 0$ (see, e.g. \cite[Chapter~1]{Knu69}), and the
properties
\begin{align*}
 \binom{n}{0} = \binom{n}{n}=1, \quad
 \binom{n}{1} = \binom{n}{n-1}=n , \quad 
 \binom{n}{k} = \binom{n-1}{k-1} + \binom{n-1}{ k},
\end{align*}
for integers $n, k \geq 1$, we obtain a structured observability
matrix
\begin{align*}
  Ob(\hat{S}, \hat{L}) 
  &= -
  \begin{bmatrix}
    \binom{0}{ 0} & \binom{1}{0} & \binom{2 }{ 0}& \cdots & \binom{r-1}{0}
    \\[0.3em] 
    \binom{1}{ 0} & \binom{2}{1} &  \binom{3 }{ 2}&\cdots & \binom{r}{r-1}
    \\[0.3em] 
    \binom{2}{ 0} & \binom{3}{1} &  \binom{4 }{ 2}&\cdots & \binom{r+1}{r-1}
    \\[0.3em] 
    \vdots & \vdots&\vdots & \ddots &\vdots \\[0.3em]
    \binom{r-1}{0} &\binom{r}{1} &\binom{r+1 }{ 2}& \cdots & \binom{2r-2}{r-1}
  \end{bmatrix} \\
  &= -
  \begin{bmatrix}
    1 & 1 & 1 &  \cdots & 1 \\
    1 & 2 & 3 & \cdots &  r \\
    1 & 3 & 6 & \cdots & \frac{r(r+1)}{2} \\
    \vdots & \vdots & \vdots & \ddots &\vdots \\
    1 &r & \frac{r(r+1)}{2}  &\cdots & \binom{2r-2 }{ r-1}
 \end{bmatrix},
\end{align*}
that is known as the Pascal matrix. It can be shown (see, e.g. \cite{BraP92}),
that the LU
decomposition of this matrix leads to its triangular factors being triangular Pascal
matrices, and thus the determinant is always 1. Consequently, $Ob(\hat{S},
\hat{L})$ always has full rank and we have established the equivalence to moment
matching choosing the expansion points as eigenvalues of $\hat{S}$,
i.e. $s_0=1$. 

\subsubsection{Jacobi Polynomials (Including Legendre and Chebychev polynomials)}
In this case we choose Sylvester equation \eqref{eq:sylv_new} to avoid
problems with a singular $\hat{E}$.
As for the Hermite polynomials we consider $\tilde{L}~=~ 
\begin{bmatrix}
 1 & 0& \hdots & 0
\end{bmatrix}.
$

For this class of orthogonal polynomials, the $\tilde{S}$ matrix has the same
structure and only differs in its entries $\alpha_i, \beta_i$ and $\gamma_i$ for
$i=1, \hdots , r$:

\begin{align*}
  \tilde{S}=
  \begin{bmatrix}
    \beta_1  & \alpha_1 &          &          &            \\
    \gamma_2 & \beta_2  & \alpha_2 &          &           \\
             & \ddots   & \ddots   & \ddots   &   \\
             &          & \gamma_{r-1}   & \beta_{r-1}   & \alpha_{r-1} \\
             &          &          & \gamma_r & \beta_r \\
  \end{bmatrix}.
\end{align*}

We now prove the full rank of the observability matrix by induction.

\textbf{Base clause:}
If we compute the observability matrix for certain small $r$ with the above
mentioned $\tilde{L}$ and $\tilde{S}$, we can clearly 
see a structure, namely:

    \begin{table}[ht]
    \centering
	\begin{tabular}{c|c|c}
	    $r$ & $Ob(\tilde{S}, \tilde{L})$ & rank \\
	    \hline 
	    \rule{0pt}{15pt}$1$ & 
	    $\begin{bmatrix}
		1
	    \end{bmatrix}$ &
	    $1 $\\
	    \hline 
	    \rule{0pt}{20pt}$2$ & 
	    $\begin{bmatrix}
		1 & 0 \\ 
		\beta_1 & \alpha_1
	    \end{bmatrix}$ & 
	    $2 $\\[1.5ex]
	    \hline 
	    \rule{0pt}{28pt}$3$ & 
	    $\begin{bmatrix}
	      1 & 0 & 0 \\
	      \beta_1 & \alpha_1 & 0 \\ 
	      \beta_1^2 + \alpha_1 \gamma_2 & \alpha_1\left(\beta_1 + \beta_2 \right) & 
	      \alpha_1 \alpha_2 
	    \end{bmatrix}$ & 
	    $3 $
	\end{tabular}
    \end{table}
These matrices are lower triangular and obviously have full rank, since
$\alpha_i \neq 0$ by definition.

\textbf{Induction hypothesis:} 
Now assume that the observability matrix of 
size $r \times r$ has a lower triangular structure 
\begin{align*}
  Ob(\tilde{S}, \tilde{L}) =
  \begin{bmatrix}
    1      &          &        &     \\
    *      & \alpha_1 &        &     \\
    \vdots &  \ddots  & \ddots &     \\
    *      & \cdots   & *      & \prod\limits_{i=1}^{r-1} \alpha_i
  \end{bmatrix}
\end{align*}
and full rank with $*$ representing the non-zero entries. 

\textbf{Induction step:} 
To prove the lower triangular structure and the full rank of the observability
matrix of size $(r~+~1) \times (r~+~1)$, we only need to have a closer look at
the $r+1$st row and column since the $r \times r$ block is unchanged.
Since the first row is the first unit vector by definition and 
$\tilde{S}$ is a tridiagonal matrix, the first $r$ entries of the last
columns are equal to zero.
The last row is computed by multiplying the $r$th row with $\tilde{S}$.
Due to the tridiagonal structure the first $r$ entries are non-zero if 
$\alpha_i, \beta_i $ and $\gamma_i$ are non-zero for $i = 1, 
\hdots , r$. The $r+1$st entry is obtained by multiplying
\begin{align*}
 \begin{bmatrix}
  * & \cdots & * & \prod\limits_{i=1}^{r-1} \alpha_i & 0
 \end{bmatrix}
 \cdot 
 \begin{bmatrix}
  0 \\ \vdots \\ 0 \\ \alpha_r \\ \beta_{r+1}
 \end{bmatrix}
 = \prod\limits_{i=1}^{r} \alpha_i .
\end{align*}
Thus, the observability matrix is a lower triangular matrix whose diagonal
entries are non-zero since $\alpha_i \neq 0 \ \forall i$, which proves its full
rank and completes the induction.

The proof for Legendre and Chebychev polynomials is omitted since these
polynomials are special cases of the Jacobi polynomials (see, e.g. 
\cite[Chapter~22]{AbrS84}).
 
\begin{figure}[t]
  \centering
  \includegraphics{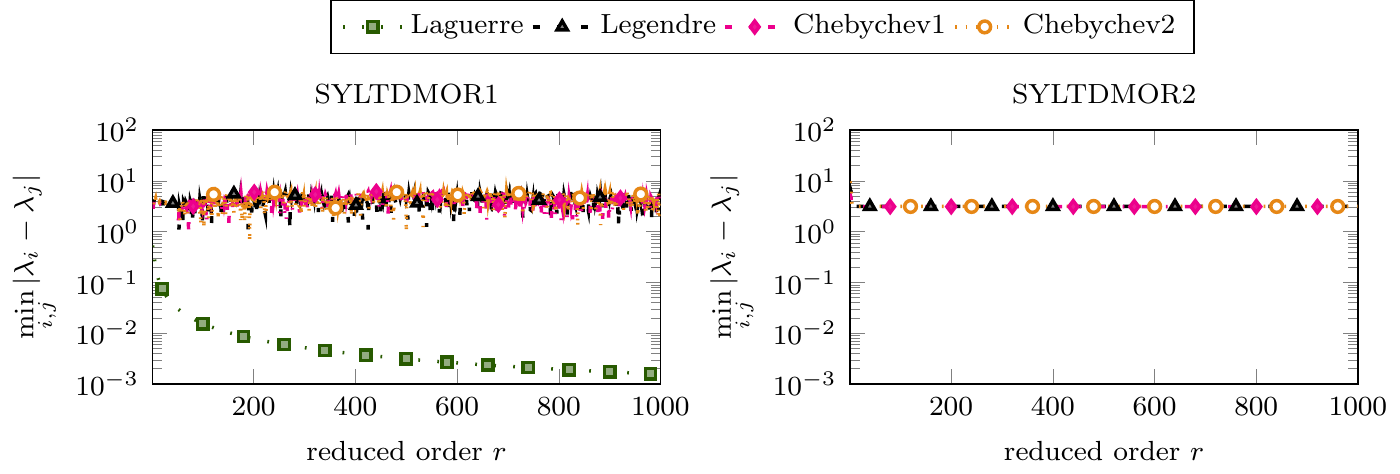}
  \caption{Minimal distance of generalized eigenvalues $(-\tilde{A},\tilde{E})$
  and $(-I_r, \hat{E})$}
  \label{fig:min_dist_eig}
\end{figure}
 
For SYLTDMOR2 presented in Section \ref{subsec:variation} we can apply 
Theorem \ref{thm:equiv_single} directly 
for the Hermite and Laguerre polynomials for SYLTDMOR2, 
since all generalized eigenvalues
of $(-I_r, \hat{E})$ are $\infty$ and $1$ respectively, and obtain an equivalence.  
Regarding the Legendre,
Chebychev polynomials of first and second kind, the Jacobi polynomials
in both approaches and the Laguerre polynomials for the original approach
proposed in \cite{morJiaC12}, we
first have to prove the distinct eigenvalues before applying Theorem
\ref{thm:equiv_multiple}. Since this is hard to verify, we show in
Figure~\ref{fig:min_dist_eig} in the left subfigure the minimal distance between
the generalized eigenvalues of $(-\tilde{A},\tilde{E})$ with matrices arising 
in Sylvester equation \eqref{eq:Sylv} in SYLTDMOR1 and  in the right subfigure for the matrix pencil
$(-I_r, \hat{E})$, where $\hat{E}$ occurs in Sylvester equation \eqref{eq:Sylv_2}
in the SYLTDMOR2 algorithm,
for a reduced order $r=1, \hdots, 1000$.  A reduced order $r~>~1000$ is
not desirable, since reducing the original system and simulating the ROM is
done with dense matrices of dimension $r$ and thus the computational costs 
become $O(r^2)$ or even $O(r^3)$, when using an implicit solving scheme,
which is way too expensive. 
Since, in the right subfigure, the minimal distance between 
the eigenvalues up to a reduced 
order of $r=1000$ is clearly larger than zero, even in finite precision, 
we conclude that one can apply Theorem \ref{thm:equiv_multiple} for the 
practically relevant $r$.  Considering the
Jacobi polynomials, we assume that $a,b >-1$ are chosen such that all
eigenvalues are distinct.
 
Summarizing these findings, we have proven the equivalence to a Krylov subspace MOR
method for SYLTDMOR2 presented in Section \ref{subsec:variation}.

\begin{theorem}[Equivalence between SYLTDMOR2 and moment matching]
\label{thm:equiv_2}
 Consider \eqref{eq:LTI} with piecewise constant input.
 
 If we choose
 \begin{itemize}
    \item Hermite polynomials,
    \item Laguerre polynomials,
    \item Legendre polynomials,
    \item Chebychev polynomials of first kind, 
    \item Chebychev polynomials of second kind,
    \item Jacobi polynomials, choosing $a$ and $b$ such that all eigenvalues of 
    $(-I_r, \hat{E})$ are distinct,
 \end{itemize}
then SYLTDMOR2 presented in Section \ref{subsec:variation} is equivalent to 
the moment matching method, where the expansion points are chosen to be the 
eigenvalues of the matrix pencil $(-I_r, \hat{E})$ arising in Sylvester equation
\eqref{eq:Sylv_2} and depend on the choice of the orthogonal polynomials and the 
reduced order.
\end{theorem}

Nevertheless, we have also shown the observability numerically in the left
subfigure of Figure~\ref{fig:rank} for the Laguerre, Legendre and
Chebychev polynomials of first and second kind. In the left subfigure of
Figure~\ref{fig:min_dist_eig} we can also see, that the minimal distance between
the generalized eigenvalues of $(- \tilde{A},\tilde{E})$ is clearly larger than
zero in case of Legendre and Chebychev polynomials of first and second kind.
Hence, for these polynomials the eigenvalues are distinct. In case of the
Laguerre polynomials the minimal distance seems to decrease, but for a reduced
order of $r=1000$ the minimal distance is approximately $10^{-3}$. 
It is thus too large for a computational error from round off accumulation. 
Hence, also the generalized eigenvalues for the Laguerre polynomials are distinct
for all relevant reduced orders \(r\). Thus, we
conjecture the equivalence to moment matching for the TDMOR approach 
of Jiang and Chen \cite{morJiaC12} and thus SYLTDMOR1:
 
 \begin{conjecture}[Equivalence between SYLTDMOR1 and moment matching]
 \label{thm:equiv_1}
  Consider \eqref{eq:LTI} with piecewise constant input and zero initial state.
 
 If we choose
 \begin{itemize}
    \item Laguerre polynomials,
    \item Legendre polynomials, odd reduced order,
    \item Chebychev polynomials of first kind, odd reduced order,
    \item Chebychev polynomials of second kind, odd reduced order,
 \end{itemize}
then SYLTDMOR1 presented in Section \ref{subsec:original} is equivalent to 
the moment matching method, where the expansion points are chosen to be the 
eigenvalues of the matrix pencil $(-\tilde{A}, \tilde{E})$ arising in Sylvester 
equation \eqref{eq:Sylv} and depend on the choice of the orthogonal polynomials
and the reduced order.
\end{conjecture}

\begin{remark}
  The equivalence between the Laguerre based time domain MOR and moment matching
  has been proven rigorously by Eid in \cite{morEid09}.
\end{remark} 

\subsection{Advantages of Moment Matching over \emph{SYLTDMOR(1/2)} in Practice}
\label{sec:advant-moment-match}

\begin{figure}[t]
  \centering
  \includegraphics{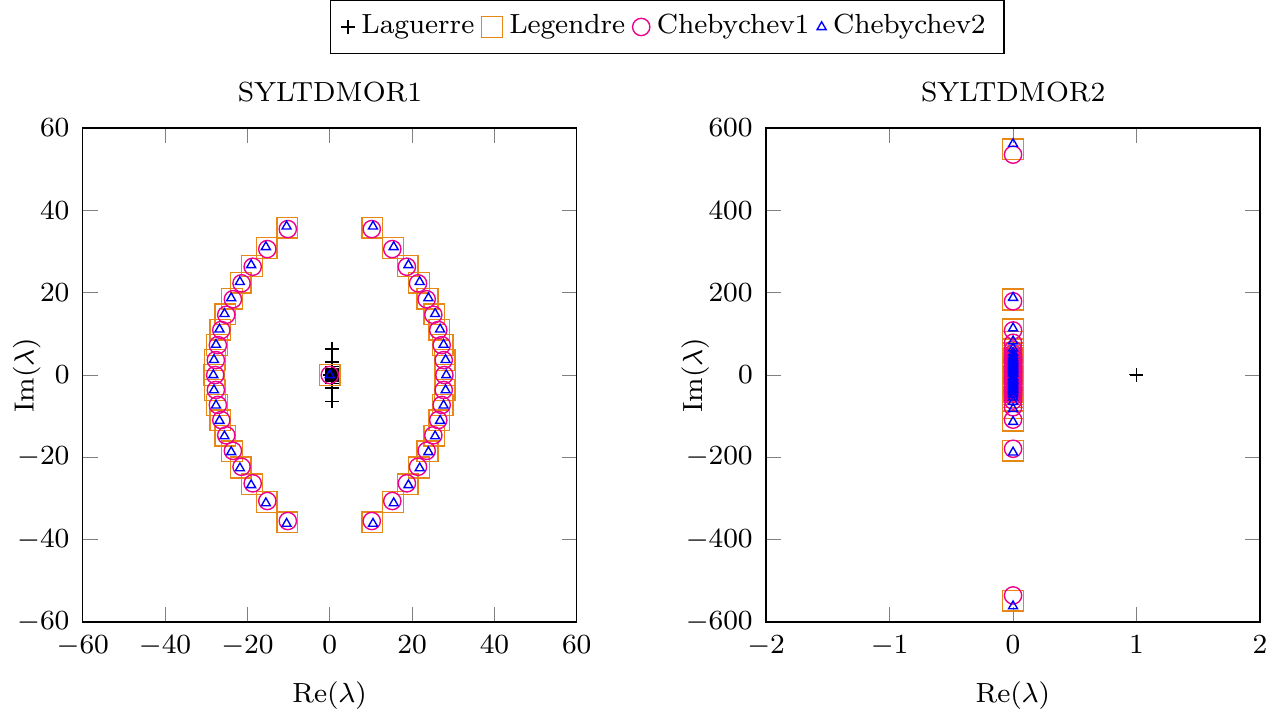}
  \caption{Generalized eigenvalues of $(-\tilde{A}, \tilde{E})$ and 
  $(-I_r, \hat{E})$ for $r=40$}
  \label{fig:eig}
\end{figure}

Comparing the polynomial based time domain MOR with moment matching one should
keep in mind, that the expansion points can be freely chosen using moment
matching. Since we use an IRKA based algorithm to determine the rational Krylov
subspace, these expansion points will be optimized fitting to the original LTI
system. Computing the reduced system with the time domain framework of Jiang
and Chen presented in Section~\ref{sec:time_domain} or SYLTDMOR2 in 
Section~\ref{subsec:variation}, the expansion points are fixed
by the choice of the family of polynomials. Exemplary, these expansion points 
or equivalently the generalized eigenvalues of matrix pairs
$(-\tilde{A}, \tilde{E})$ (SYLTDMOR1) and $(-I_r, \hat{E})$ 
(SYLTDMOR2), are shown in Figure~\ref{fig:eig} for a reduced order $r=40$. 

Thus, a variation of the expansion points arising from the time
domain approach can be achieved by in- or decreasing the reduced order $r$ or
rather using another family of orthogonal polynomials, on the one hand.
On the other hand, orthogonal polynomials of higher degree $s \in \mathbb{Z}$
up to degree $s+r-1$ could be used or $r$ arbitrarily chosen orthogonal
polynomials. 
To the best of our knowledge this has not been tried in the open literature so 
far. 

Considering the choice of expansion points, the time domain MOR framework
based on orthogonal polynomials, thus only seems to be a
restriction of moment matching to a rather limited set of possible 
combinations.
Thus, in practice IRKA has to be expected to provide more accurate ROMs in 
almost all cases.
Further, IRKA has to be at least as good as the discussed time domain MOR 
approaches.

\section{Numerical Examples} \label{sec:example}
To illustrate the effectiveness of the time domain MOR techniques
based on orthogonal polynomials, but also to confirm our conjecture about
the restriction of this method compared to moment
matching, we will present results for three well-known test examples.

The basic setup for the first two examples is the same, namely:
\begin{itemize}
    \item initial state $x_0 = \mathbb{O}_{n,1} \in \mathbb{R}^n,$
    \item time interval $t\in [0,1]$, i.e. $t_0=0$,
    \item time step $\tau= 0.001$,
    \item input $u(t)=
      \begin{cases}
	0 \qquad &,t\in [0,0.1)\\
	\frac{1}{2} \sin \left( \pi (10t- \frac{3}{2})\right) +\frac{1}{2} \quad &,
	t\in [0.1,0.2) \\
	1 &,t \in [0.2,1]
      \end{cases}.
      $
\end{itemize}
The test system is always the following:
\begin{itemize}
  \item CPU: 2x Intel\textsuperscript{\textregistered} Xeon\textsuperscript{\textregistered} X5650
  \begin{itemize}
    \item $6$ Cores per CPU,
    \item clock rate: $2.67$ GHz,
    \item $12$ MB Cache per CPU,
  \end{itemize}
\item memory: $48$ GB DDR3 with ECC.
\end{itemize}

All examples are computed using \matlab{} R2012b.

In this paper, we illustrate the 2-norm averaged relative error over time,
i.e. 
\begin{align*}
  \left\Vert \frac{y(t) - y_r(t)}{y(t)} \right\Vert_{2, [0,1]} 
  & = 
  \left(\int\limits_{0}^1 \left( \frac{y(t) - y_r(t)}{y(t)} \right)^2 dt 
  \right)^{\frac{1}{2}} \\
  & \approx 
    \left(\sum\limits_{i=1}^{\frac{1}{\tau}} 
	\left( 
	      \frac{y(i\tau) - y_r(i\tau)}{y(i\tau)}
	\right)^2 
    \right)^{\frac{1}{2}},
\end{align*}
where $y(t)$ and $y_r(t)$ are computed using the implicit Euler method (see, e.g.
\cite[Chapter~2]{But03}). Note, that the $\infty$-norm averaged error over 
time looks comparable and is thus not illustrated.
We also show the total time, that was spent to reduce the original LTI system
and to solve the reduced LTI system compared to the time, that was spent on 
solving the original LTI system. 
In these figures, we compare SYLTDMOR1 and SYLTDMOR2 to the two most 
important and well accepted methods for stable LTI systems, IRKA 
(one- and two-sided) and balanced truncation.
For the one-sided IRKA approach, the projection matrix $V$ is computed from
the output Krylov subspace, i.e. in case of one expansion point 
$s_0$ $V$ is the basis of 
$\mathcal{K}_r \left( (A-s_0E)^{-T}E^T,(A-s_0E)^{-T}C^T \right)$.
In case of multiple expansion points $s_1, \hdots , s_r $ $V$ is a basis of
$ \bigcup\limits_{i=1}^r \mathcal{K}_{r_i} \left( (A-s_iE)^{-T}E^T,
  (A-s_iE)^{-T}C^T \right)$.
This is implemented according to the theory in Section 
\ref{sec:moment_matching} and does not present a restriction to moment
matching since upon convergence of IRKA, the expansion points are (locally)
optimally placed for the system with respect to \(\mathcal{H}_{2}\) 
approximation.  In the figures, we use the following notations for Chebychev
polynomials of first (Chebychev1) and second kind (Chebychev2), the one-sided
IRKA resp. moment matching (oIRKA/oMM), the two-sided IRKA resp. moment matching
(IRKA/tMM) and balanced truncation (BT).
Note, that we computed $50$ cycles to average the results.


\subsection{Triple Peak Example}\label{subsec:ex_1}
Our first example is the triple peak, sometimes also called FOM, example (see, e.g.
\cite{morPen99}), where the dynamical system \eqref{eq:LTI} of order $n=1\,006$ 
is given by 
\begin{figure}[t]
  \centering
  \includegraphics{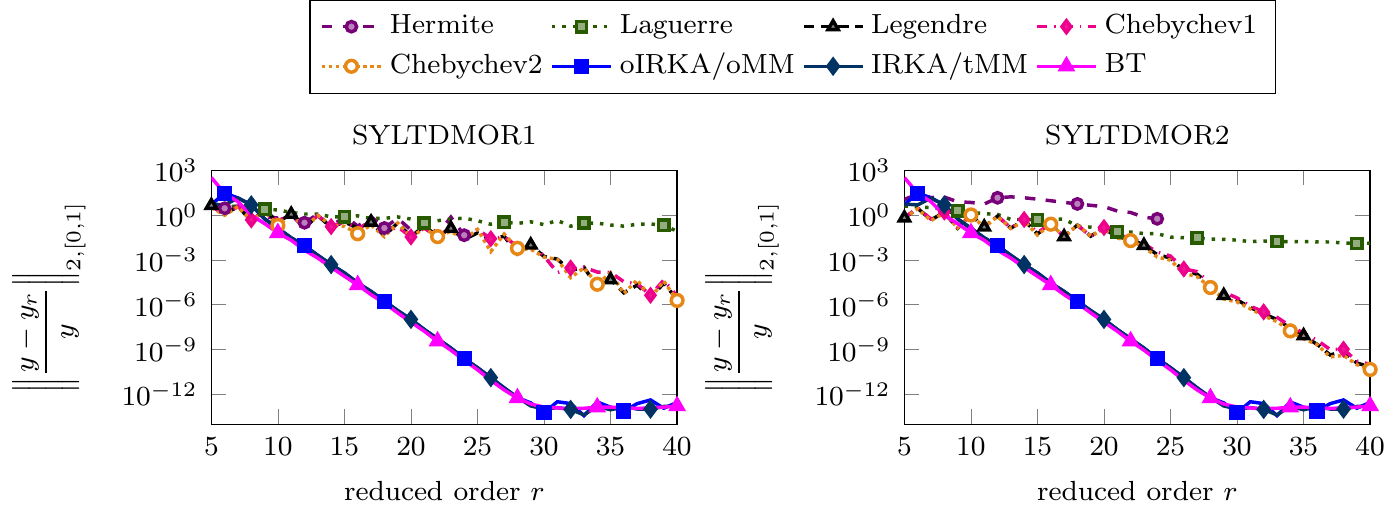}
  \caption{Relative errors (triple peak example of Section \ref{subsec:ex_1})}
  \label{fig:rel_err_peak}
\end{figure}
$E=I_n, A \in \mathbb{R}^{1 \,006 \times 1\, 006}$ a block diagonal matrix of the form
\begin{alignat*}{3}
  A & = 
  \begin{bmatrix}
    A_1 &     &     &     \\
        & A_2 &     &     \\
	&     & A_3 &     \\
	&     &     & A_4 \\
  \end{bmatrix}, \quad 
  A_1 &=
  \begin{bmatrix}
    -1 & 100 \\ -100 & -1
  \end{bmatrix},\\
  A_2 &=
  \begin{bmatrix}
    -1 & 200 \\ -200 & -1
  \end{bmatrix}, \qquad \qquad \ 
  A_3 &=
  \begin{bmatrix}
    -1 & 400 \\ -400 & -1
  \end{bmatrix},
  \\
  A_4 &= \text{diag}\{ -1, -2 , \hdots, -1\,000 \}, &&
\end{alignat*}
and the input and output matrices are 
\begin{align*}
  \nonumber
  B^T = C= [\underbrace{10, \hdots, 10}_{6}, \underbrace{1, \hdots ,1}_{1\,000}] 
        \in \mathbb{R}^{1 \times 1\,006}.
\end{align*}
That means $x(t)\in \mathbb{R}^{1 \,006}$ and $u(t), y(t) \in \mathbb{R}$.

In Figure~\ref{fig:rel_err_peak} the $2$-norm averaged relative error over time is
illustrated for the reformulated time domain approach SYLTDMOR1, 
presented in Section~\ref{subsec:original}, and its variation SYLTDMOR2, 
presented in Section
\ref{subsec:variation}, for Hermite, Laguerre, Legendre and Chebychev polynomials
of first and second kind, each compared to the one- and two-sided IRKA and
balanced truncation. In both subfigures we can easily see, that the frequency 
domain MOR
approaches approximate the original system, for a reduced order $r \geq 11$,
much better than the time domain approaches. Only the Legendre and the two
Chebychev polynomial families in SYLTDMOR2 show a considerable
decay of the
relative error ending up with a relative error of $10^{-10}$ for a reduced order
$r=40$.  The same orthogonal polynomials in SYLTDMOR1 have only a
slight decay of the relative error ending up at around $10^{-5}$ for a reduced
order $r=40$. 
Compared to this, balanced truncation and both IRKA approaches show a nicer 
decay ending with a relative error of $10^{-12}$ for $r \geq 28$ and are thus
at least $2$ orders lower compared to the Legendre and both Chebychev 
polynomials for both represented time domain approaches. 
Regarding the Laguerre polynomials the relative error does not
seem to change for an increasing reduced order $r$ and stagnates around $10^{-1}$.
This phenomenon can be explained, if we look at their differential recurrence
coefficients stated in Table~\ref{tab:rec}. These coefficients are constant and do
not depend on the reduced order $r$ such that matrix $\hat{S}$ is always a 
triangular matrix with only $1$ on the diagonal as seen in Section
\ref{subsec:Laguerre}. Thus the only expansion point of these
polynomials for SYLTDMOR2, given by the multiple eigenvalue of $\hat{S}$, is $1$ with multiplicity $r$. 
That means in the sense of moment matching, that moments up 
to order $r$ are matched. But on the other hand all other important frequencies 
are ignored leading to a bad approximation. Considering SYLTDMOR1 it might be possible
that this effect is also caused by the clustering of the eigenvalues.
A special case are the Hermite
polynomials, since they are only illustrated until $r=24$. The reason for this
is the extreme condition number, i.e. numerical singularity, of the matrices $H$
and $\hat{H}$, which makes it impossible to solve the linear system of
equations. 
\begin{figure}[t]
  \centering
  \includegraphics{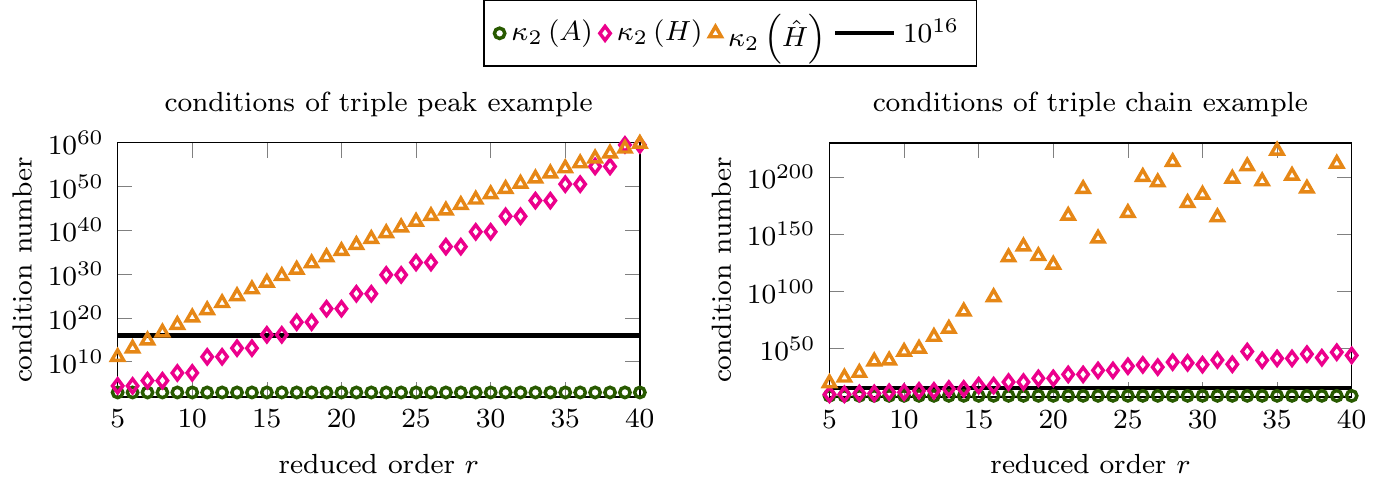}
  \caption{Condition numbers using Hermite polynomials}
  \label{fig:cond}
\end{figure}
To get an impression how the condition numbers grow by increasing
the reduced order, the $2$-norm condition numbers $\kappa_2(.)$ using Hermite 
polynomials are shown for the triple peak example presented in Section
\ref{subsec:ex_1} and 
the example of the following Section \ref{subsec:ex_2} for the matrices $A, H$ 
and $\hat{H}$ in Figure~\ref{fig:cond}. Since in \matlab{} a matrix is 
numerically not invertible for $\kappa_2(.)> 10^{16}$, this bound is added in 
the subfigures of Figure~\ref{fig:cond}.
Nevertheless, we can easily see in Figure~\ref{fig:rel_err_peak}, that
the Hermite polynomials have a large relative error of around $10^{-1}$ in 
SYLTDMOR1 and about $10^0$ in SYLTDMOR2 at
$r=24$, such that the original system is not approximated well.

\begin{figure}[!t]
  \centering
  \includegraphics{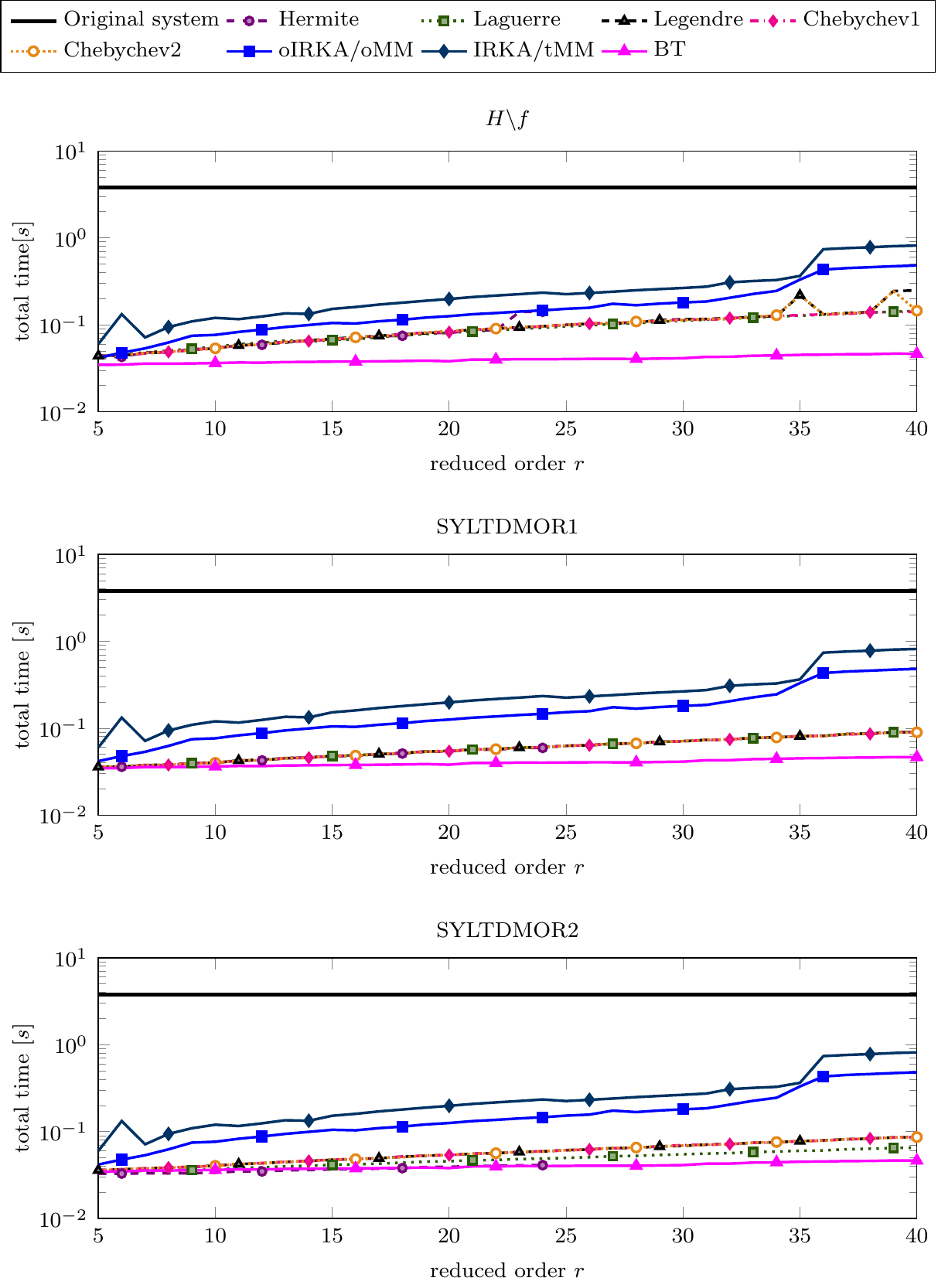}
  \caption{Total time (triple peak example of Section \ref{subsec:ex_1})}
  \label{fig:time_peak}
\end{figure}

The total time, that is spent on reducing the original system 
and solving the reduced system, is
illustrated in Figure~\ref{fig:time_peak} and is compared to the time, that is
needed to solve the original $(1\,006 \times 1\,006)$ system. This figure is
divided into three subfigures containing the reduction and solution times for
determining $V$ by solving the huge $(nr \times nr)$ linear system of equations
\eqref{eq:GLS} with \matlab{}'s backslash operator and Sylvester
equations \eqref{eq:Sylv} and \eqref{eq:Sylv_2} with the method from
\cite{morBenKS11}. In all three subfigures of Figure~\ref{fig:time_peak} it is
easy to see, that all solution methods are faster than solving the original
system. Comparing the time with the backslash operator and solving a Sylvester
equation, irrespective whether SYLTDMOR1 or SYLTDMOR2 is chosen,
the Sylvester solver is, for $r=40$, up to two times faster than solving with
\matlab{}'s backslash.  For a reduced order of $r \geq 34$ there is an
increase of the total time for both IRKA variants. The reason for this can be
found in Figure~\ref{fig:rel_err_peak}, since the relative error of these methods
is close to machine precision and thus an improvement of the relative error is
not possible.  

Basically, this example demonstrates, that the time domain MOR 
framework~\cite{morJiaC12} reduces the LTI system fast.
Further, these time domain MOR approach might be even faster by inserting 
the eigenvalues of matrix pencils $(-\tilde{A},\tilde{E})$ or $(-I_r, \hat{E})$ 
directly as expansion points in the moment matching method ending up with 
the same projection matrix according to Theorem \ref{thm:equiv_2} and Conjecture
\ref{thm:equiv_1}. This could be implemented by precomputing the generalized 
eigenvalues and saving them in a data base, such that they are quickly available.
Here, the Jacobi polynomials are excluded, since they require a further analysis
of the choice of parameters $a$ and $b$.
Nevertheless, this example also demonstrates, that the time domain MOR approaches
are less accurate compared to IRKA or balanced truncation.

 
\subsection{Triple Chain Example}\label{subsec:ex_2}
The second example is the triple chain example from \cite{TruV09},
i.e. three mass-spring-damper chains of length $200$ are fixed by one coupling
mass. Since this example, which is parametrized as in \cite{Saa09}, results in a
second order systems
\begin{align}
\label{eq:sec_ord}
    \begin{split}
      M\ddot{x}(t) + D \dot{x}(t) + K x(t) &= \mathcal{B} u(t), \\
      y(t)&= \mathcal{C} x(t),
    \end{split}
\end{align}
we transform it into the first order system 
\begin{align*}
  \underbrace{\begin{bmatrix} K & 0 \\ 0 & M\end{bmatrix}}_{E}
  \underbrace{\begin{bmatrix} \dot{x}(t) \\ \ddot{x}(t)\end{bmatrix}}_{\dot{z}(t)}
  &= \underbrace{\begin{bmatrix}0 & K \\ -K & -D\end{bmatrix}}_{A}
  \underbrace{\begin{bmatrix} x(t) \\ \dot{x}(t)\end{bmatrix}}_{z(t)}
  +
  \underbrace{\begin{bmatrix} 0 \\ \mathcal{B}\end{bmatrix}}_{B}
  u(t),\\
  y(t) &=
  \underbrace{\begin{bmatrix}\mathcal{C} & 0\end{bmatrix}}_{C}
  \underbrace{\begin{bmatrix} x(t) \\ \dot{x}(t)\end{bmatrix}}_{z(t)}.
\end{align*}
Here, the matrices of the second order system \eqref{eq:sec_ord} are given by a
diagonal matrix $M \in \mathbb{R}^{601 \times 601}$ containing the masses, the
damper matrix $D \in \mathbb{R}^{601 \times 601}$ and the stiffness matrix
$K \in \mathbb{R}^{601 \times 601}$. The input and output matrices are again
transposes of each other given by
$\mathcal{B}^T= \mathcal{C}=[1 \hdots 1] \in \mathbb{R}^{1 \times 601}$.  Hence,
$E, A \in \mathbb{R}^{1\,202 \times 1\,202},B^T,C \in \mathbb{R}^{1 \times 1\,202},
z(t)\in \mathbb{R}^{1\,202}$ and $ u(t),y(t) \in \mathbb{R}$.
 
\begin{figure}[t]
  \centering
  \includegraphics{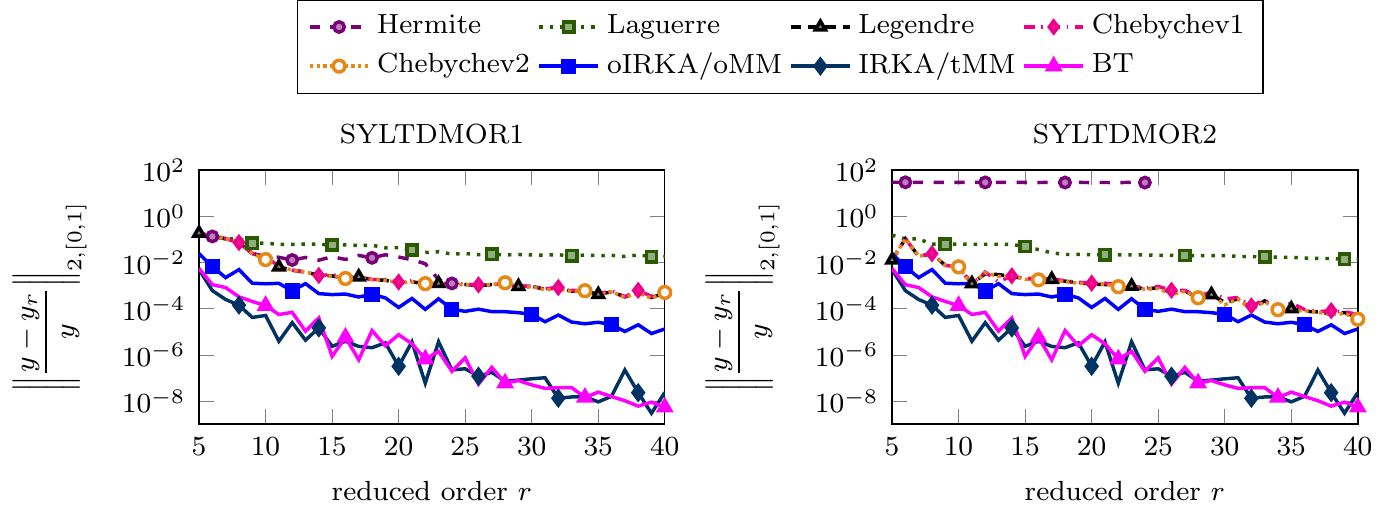}
  \caption{Relative errors (triple chain example of Section \ref{subsec:ex_2})}
  \label{fig:rel_err_chain}
\end{figure}

Figure~\ref{fig:rel_err_chain} illustrates the $2$-norm averaged relative error
over time for the time domain MOR approaches compared to the frequency
domain MOR methods
as mentioned in the previous example. Again, the relative error using Laguerre
polynomials shows only minor changes and is for a reduced order $r=40$ at around
$10^{-1}$ using SYLTDMOR1 or SYLTDMOR2. Similar to the Laguerre
polynomials, the Hermite polynomials do not approximate the original system well
since the relative error is, especially in case of SYLTDMOR2, too large, 
namely $29$ for a reduced order $r = 24$. In case of
SYLTDMOR1 the relative error decreases to around $10^{-2}$, but since
the $H$ matrix is numerically singular, it is not possible to determine
further projection matrices.  In both subfigures the 2-norm averaged relative
error over time decreases, if Legendre or Chebychev polynomials of first or
second kind are used. This relative error is around $10^{-3}$ for SYLTDMOR1
and $10^{-4}$ for SYLTDMOR2, for a reduced order $r=40$.  Comparing
both time domain approaches to the one-sided IRKA algorithm, we can clearly see,
that the relative error using moment matching is always at least as small as in
case of the time domain MOR, but mostly even smaller. Looking at the 
two-sided IRKA algorithm and the balanced truncation method, we see, 
that these methods are even more successful since they have a steeper 
decrease of the relative error ending up with a relative error of
order $10^{-8}$, which is $4$ orders of magnitude lower 
than for the one-sided method.
 
 \begin{figure}[!t]
  \centering
  \includegraphics{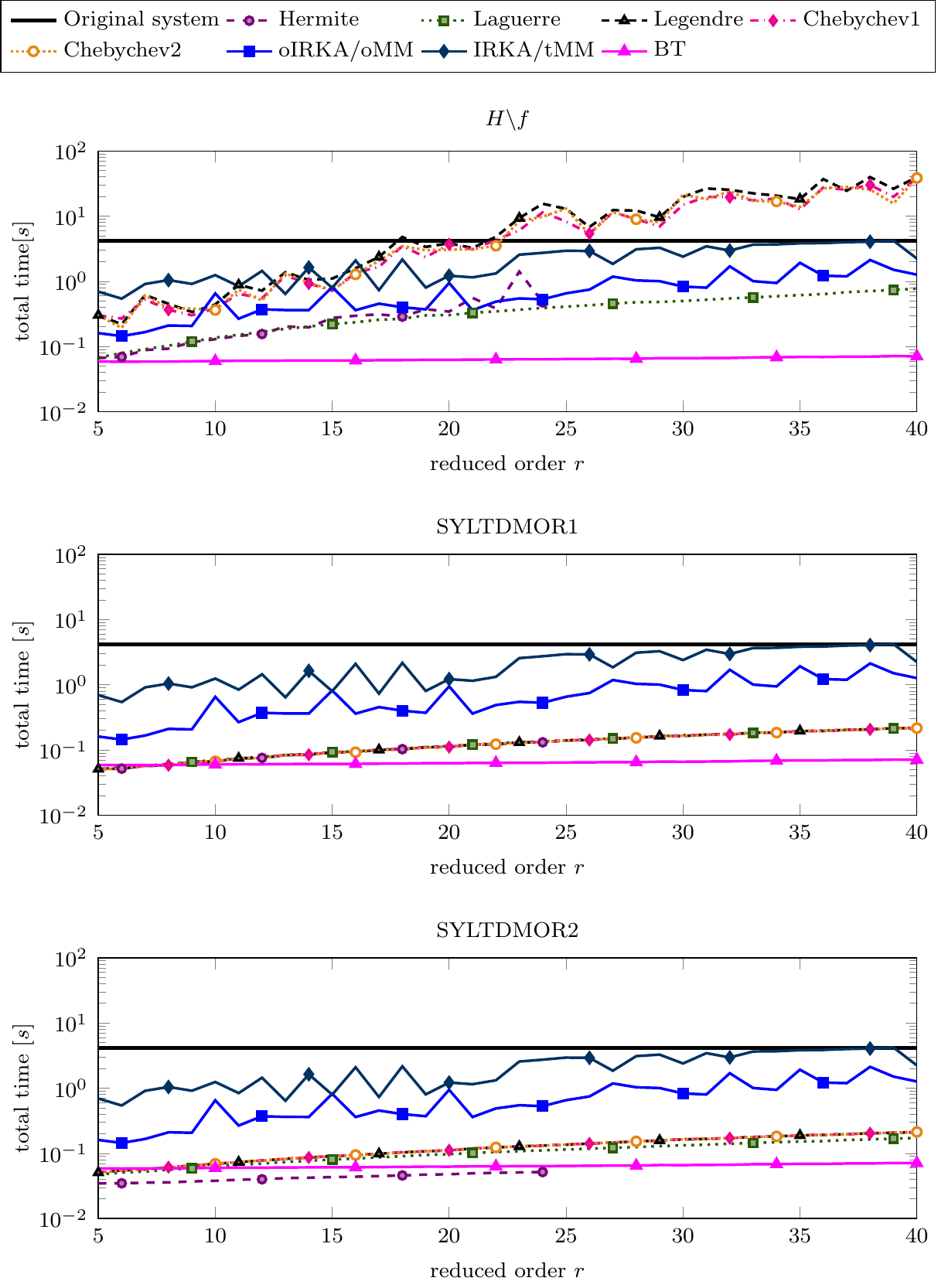}
  \caption{Total time (triple chain example of Section \ref{subsec:ex_2})}
  \label{fig:time_chain}
\end{figure}

We now take look at the total time, i.e. the time needed to reduce the 
original system and to simulate the reduced system.
In the upper subfigure of Figure~\ref{fig:time_chain}, illustrating the total 
time for the backslash solver, we can easily see that, if we use the time 
domain approach with Hermite and Laguerre polynomials, we are faster than 
solving the original system.
Unfortunately, these polynomials do not approximate the original system well. 
Using the Legendre or one of the Chebychev families, this approach is 
faster than simulating the original system until $r=21$. 
Hence, considering a larger reduced order is not effective any more. 
If we use one of the Sylvester solvers instead, we can achieve a high
speed-up, e.g. choosing a reduced order $r=40$, we can reduce the system and
simulate the reduced system up to $19$ times faster than simulating
the original system. 
This also holds for the Hermite and Laguerre polynomials. 
The former ones are only considered up to a reduced order $r=24$.
Regarding both IRKA implementations and balanced truncation, it is clearly
visible, that these MOR techniques consume more time than the time domain
approaches using Sylvester equations. 
Nevertheless, these methods are still faster than simulating the original triple 
chain example.


\subsection{Butterfly Gyroscope Example}\label{subsec:ex_3}

A more practice oriented and larger example is given by the butterfly gyroscope
example from the Oberwolfach benchmark collection for model order reduction
(see, \cite{morwiki_gyro}) described in \cite{morBil05}, which represents a vibrating
micro-mechanical gyroscope. The device consists of a three-layer silicon
wafer stack. Its middle layer contains the sensor element, which consists of two
wing pairs that are connected to a common frame -- the reason the gyro is called
butterfly. This example is given as second order system \eqref{eq:sec_ord} and
transformed into a first order realization
\begin{align*}
      \underbrace{
      \begin{bmatrix}
	- K^T & 0  \\ 0 & M
      \end{bmatrix}
      }_{E}
      \underbrace{
      \begin{bmatrix}
	\dot{x}(t)  \\ \ddot{x}(t)
      \end{bmatrix}
      }_{\dot z(t)}
      & = 
      \underbrace{
      \begin{bmatrix}
	0 & - K^T  \\ -K & -D
      \end{bmatrix}
      }_{A}
      \underbrace{
      \begin{bmatrix}
	x(t)  \\ \dot{x}(t)
      \end{bmatrix}
      }_{z(t)}
      +
      \underbrace{
      \begin{bmatrix}
	 0  \\ \mathcal{B}
      \end{bmatrix}
      }_{B}
      u(t) \\
      y(t) 
      & =
      \underbrace{
      \begin{bmatrix}
	\mathcal{C} & 0 
      \end{bmatrix}
      }_{C}
      \underbrace{
      \begin{bmatrix}
	x(t)  \\ \dot{x}(t)
      \end{bmatrix}
      }_{z(t)}.
\end{align*}
Here, matrices $M,D,K \in \mathbb{R}^{17\,361 \times 17\,361}$ of the second order system 
are symmetric matrices, the input is the vector $\mathcal{B} \in
\mathbb{R}^{17\,361}$ and output matrix 
is given as $\mathcal{C} \in \mathbb{R}^{12 \times 17\,361}$. In order to obtain a SISO system, 
we only consider the first row of $\mathcal{C} $.
Hence the matrices and vectors of the first order system are of dimension
$E ,A \in \mathbb{R}^{34\,722 \times 34\,722} ,
B^T, C \in \mathbb{R}^{1 \times 34\,722}, z(t) \in \mathbb{R}^{34\,722}$ and $u(t),
y(t) \in \mathbb{R}$.

Accounting for the output trajectory, we changed the computational setup for this
example to 
\begin{itemize}
 \item time interval $t \in [0,0.005]$,
 \item time step $\tau = 5 \cdot 10^{-6}$, and
 \item input \(u(t)\) the corresponding sine-smoothed step from 0 to 1 in the
   interval \([0.0005, 0.001]\). 
\end{itemize}
However, this can easily be transformed to a time interval $t \in [0,1]$ using the transformation
$t_1 = \frac{t_1}{0.005}$ with time step $\tau = 0.001$ as in the previous examples.
The initial state $x_0$ and the input $u(t)$
are chosen as in the previous examples. Furthermore, we computed only one cycle instead
of $50$ to limit the computation time, since external effects are expected to be
far less influential in this case. 

\begin{figure}[!t]
  \centering
  \includegraphics{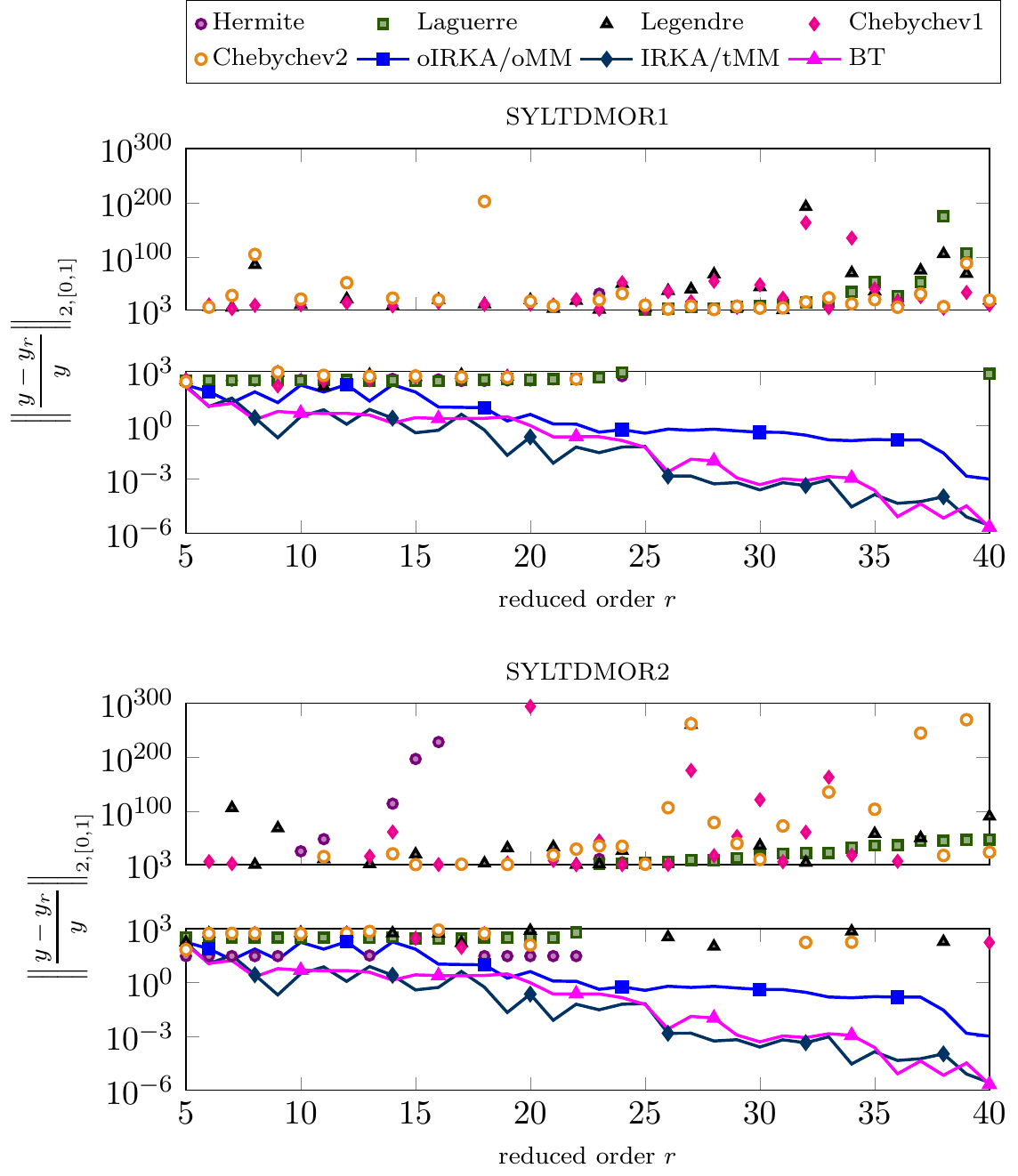}
  \caption{Relative errors (butterfly gyroscope example of Section \ref{subsec:ex_3})}
  \label{fig:rel_error_gyro_splitted_axis}
\end{figure}

In Figure \ref{fig:rel_error_gyro_splitted_axis} the $2$-norm averaged relative error
over time is illustrated. These subfigures differ from the ones from the previous examples.
Here, each subfigure consists of two plots, where the upper plot is logarithmically scaled
from $10^3$ to $10^{300}$ and the lower plot is logarithmically scaled from $10^{-6}$ to 
$10^3$. 
The reason for this unusual scaling are large differences of the relative errors of the 
time domain MOR approach using orthogonal polynomials compared to the remaining MOR methods 
as IRKA and balanced truncation.
While the maximum relative error of the remaining methods is around $10^3$, the minimum 
relative error of the orthogonal polynomials is in the same area and increases even up to
an order of $10^{294}$. 
Furthermore it is even possible, that this method produces unsuitable, i.e. non
stable ROMs. Then the relative error is given by \texttt{NaN} and omitted in the
plot.
Compared to this, the relative error for the balanced truncation method and both IRKA 
approaches nicely decreases. Looking at reduced order $r=40$, balanced truncation and the two-sided
IRKA method end up with a relative error of around $10^{-5}$ and the one-sided IRKA approach
with a relative error of around $10^{-2}$.
Thus Figure~\ref{fig:rel_error_gyro_splitted_axis} illustrates impressively, that the 
orthogonal polynomial based time domain MOR framework fails completely for this example.

\begin{figure}[t]
  \centering
  \includegraphics{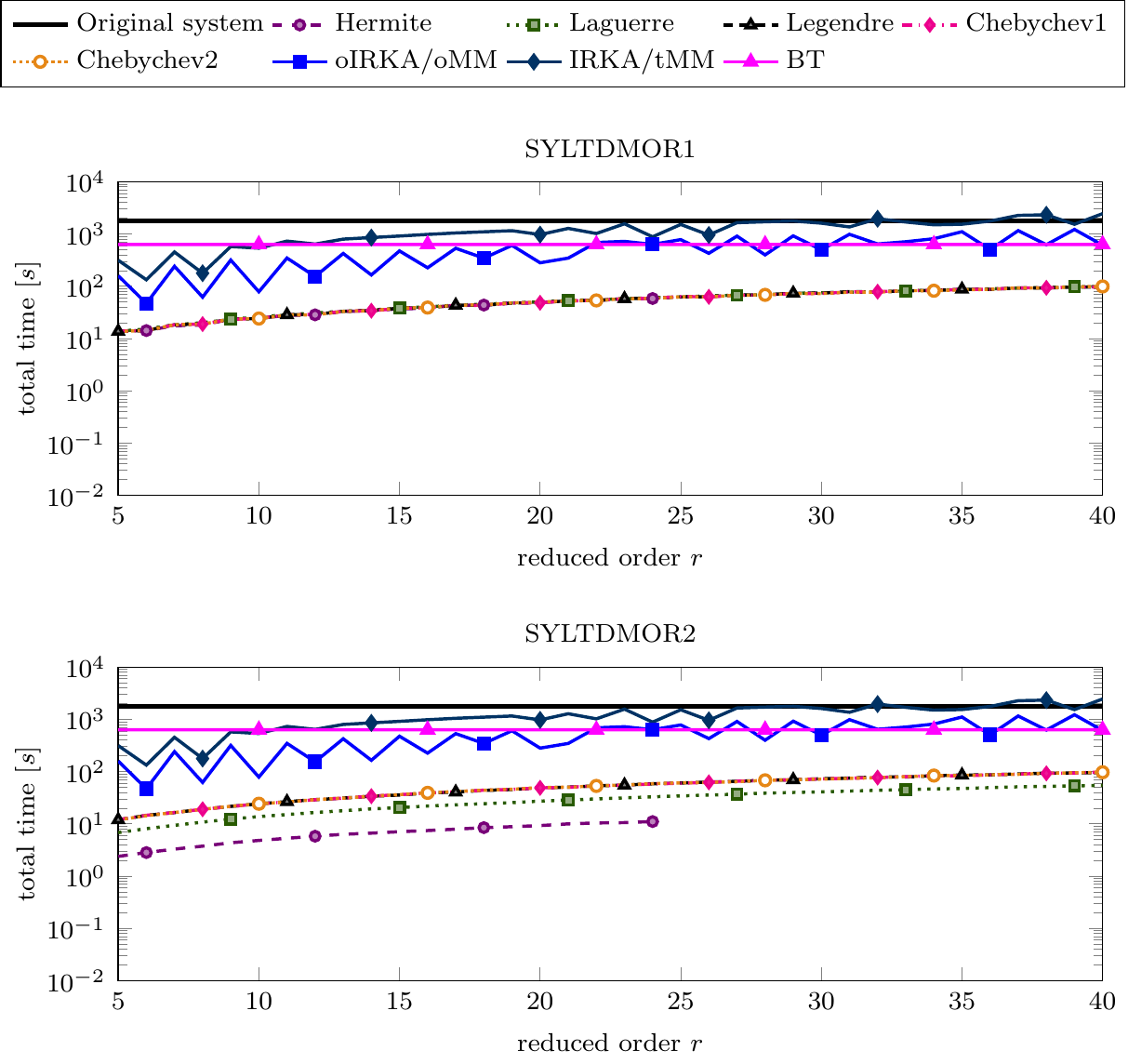}
  \caption{Total time (butterfly gyroscope example of Section \ref{subsec:ex_3})}
  \label{fig:time_gyro}
\end{figure}

Figure \ref{fig:time_gyro} illustrates the total time that is needed to reduce and 
simulate the reduced system. Just like in the previous examples, SYLTDMOR1 and SYLTDMOR2 are clearly
faster than using one of the IRKA approaches or the balanced truncation
method. Still, the increased computation time is a price worth paying, since the
latter methods produce suitable and reliable ROMs in all test cases. 
Note, that we cannot report the total time using \matlab{}'s backslash
operator, since the computations became too memory consuming even for our well
equipped test system.

If we now take the relative errors illustrated in 
Figure~\ref{fig:rel_error_gyro_splitted_axis} into account, IRKA and balanced truncation 
are clearly more desirable than the time domain MOR techniques, since all of them
approximate the original model behavior far better.

\begin{figure}[t]
  \centering
  \includegraphics{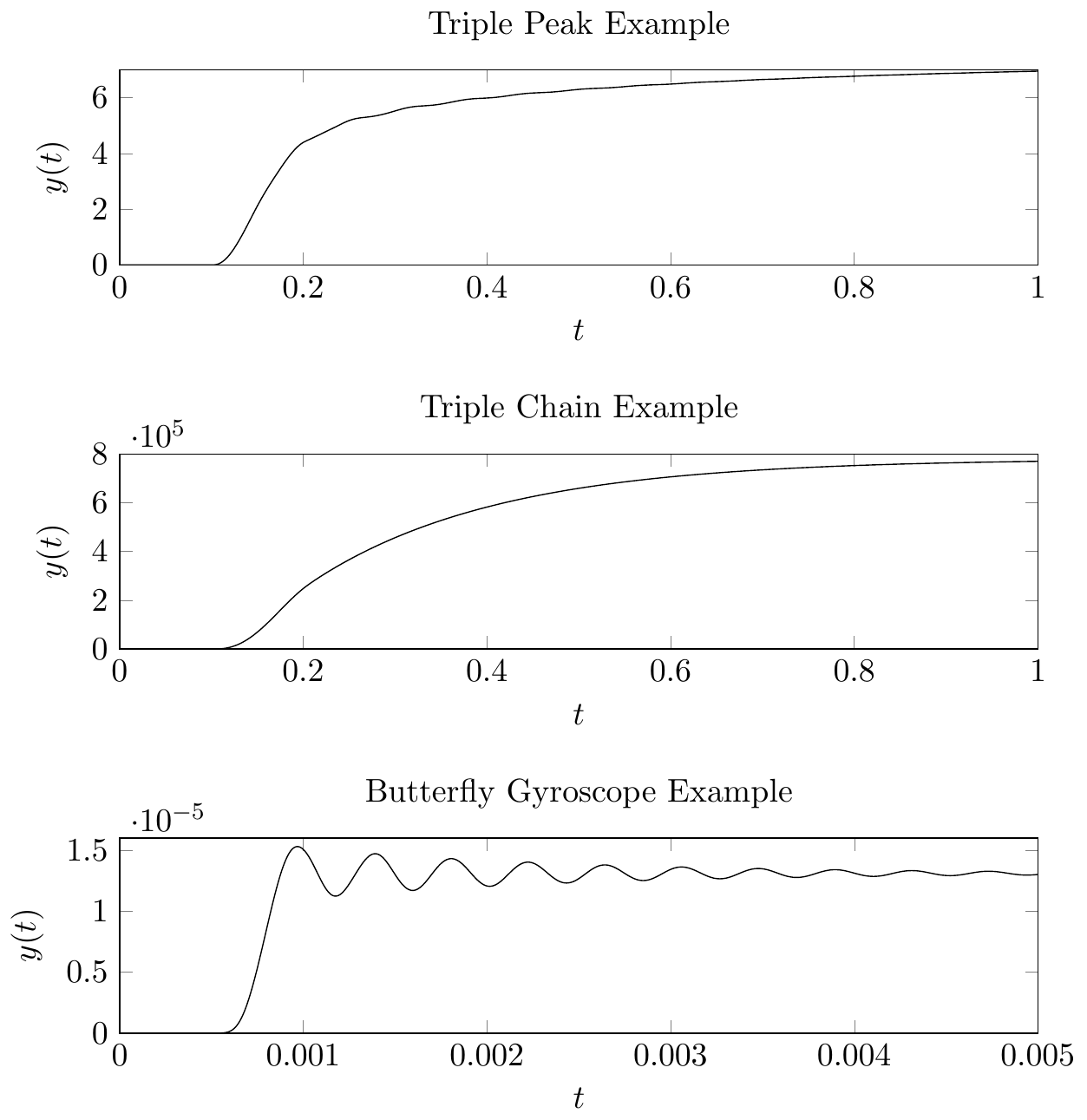}
  \caption{Smoothed step-response trajectories of the three test examples.}
  \label{fig:trajectories}
\end{figure}

\subsection{Observations}
The numerical examples, reported above, give us an impression of the
effectiveness of the time domain framework of Jiang and Chen and its variation
SYLTDMOR2 compared to other frequency domain MOR techniques. 
Regarding the relative errors, the one- and
two-sided IRKA algorithm and balanced truncation always have a stronger decrease
and thus (almost) always smaller relative errors. The fact, that these methods
consume more time than the time domain approaches is negligible, since these
methods (nearly) always performed faster than solving the original system.

Note, that we also tried to implement the iterative splitting method proposed in 
\cite{morJiaC12}, but our implementation never converged for the above 
mentioned examples. Thus, we are not able to compare with our proposed 
solving techniques. 
Since solving with \matlab{}s backslash operator is time-consuming,
we also tried to solve \eqref{eq:GLS} using the preconditioned Generalized Minimum
Residual (GMRES) method.
But since the matrix $H$ in \eqref{eq:GLS} is ill-conditioned, we were not able
to find a good preconditioner, such that \eqref{eq:GLS} can be solved fast.

The aim of model reduction in our context is to find a reduced system
approximating the original system well, such that the time, that is spent on
reducing and solving the reduced system, is less than simulating the original
system. SYLTDMOR1 and its variation SYLTDMOR2 clearly consume less
time than solving the original system, but the relative error either decreases very
slowly as in the first two examples and thus the reduced order needs to be comparably 
large or the computed ROM is not suitable at all. 
This behavior can be easily explained when looking at the trajectories of all examples
illustrated in Figure \ref{fig:trajectories}. The trajectories of the triple peak 
and the triple chain example can easily be expressed by using low-order
polynomials. In contrast to this, the trajectory of the butterfly gyroscope
example is oscillating rather fast and thus it requires higher order polynomial to
approximate the solution properly, which are not present in the bases generating
low-order models following the Jiang/Chen framework. 

\section{Concluding Remarks} 
\label{sec:conclusion}
In this paper, we picked up the time domain MOR framework based on the idea
of Jiang and Chen in \cite{morJiaC12} and transformed the resulting huge
linear system of equations into a Sylvester equation, that can be solved 
very efficiently. A slight variation of the formulation leads to another
even nicer Sylvester equation considering a fixed initial state and leading to
easier structure in the small coefficient matrices. 
Using the duality theorem, we show a connection of the time domain
MOR methods to moment matching, but we also illustrate, that the expansion 
points created by the time domain approaches cannot adapt to the system 
and thus the time domain approaches in this paper cannot 
keep up with proper moment matching, which is only one iteration
step of IRKA.

\section*{Code Availability}
The \matlab{} implementation used to compute the presented results
can be obtained from
\begin{center}
    \framebox{
      doi: 10.5281/zenodo.1243090
    }
\end{center}

\noindent and is authored by: {\scshape Manuela Hund} and {\scshape Jens Saak}.

\end{document}